\documentclass[12pt]{amsart}
\usepackage{hyperref}
\begin{document}
\title{A Hopf Index Theorem for foliations}
\author{Victor Belfi}
\address{\hskip-\parindent
        Victor Belfi\\
        Department of Mathematics\\
        Texas Christian University\\
        Box 298900\\
        Fort Worth, Texas 76129}
\email{v.belfi@tcu.edu}
\author{Efton Park}
\address{\hskip-\parindent
        Efton Park\\
        Department of Mathematics\\
        Texas Christian University\\
        Box 298900\\
        Fort Worth, Texas 76129}
\email{e.park@tcu.edu}
\author{Ken Richardson}
\address{\hskip-\parindent
        Ken Richardson\\
        Department of Mathematics\\
        Texas Christian University\\
        Box 298900\\
        Fort Worth, Texas 76129}
\email{k.richardson@tcu.edu, krichar@msri.org}

\date{April, 2001}

\thanks{Research at MSRI is supported in part by NSF grant DMS-9701755.}

\begin{abstract}
We formulate and prove an analog of the Hopf Index Theorem for Riemannian
foliations. We compute the basic Euler characteristic of a closed Riemannian
manifold as a sum of indices of a non-degenerate basic vector field at
critical leaf closures. The primary tool used to establish this result is an
adaptation to foliations of the Witten deformation method.
\end{abstract}

\maketitle





\renewenvironment{Sb}{_\bgroup \subarray{c}}{\endsubarray\egroup } 

\renewcommand{\frak}{\mathfrak}

\renewcommand{\Bbb}{\mathbb}

\newcommand{\func}{\mathrm}
\newtheorem{theorem}{Theorem}[section]
\newtheorem{corollary}[theorem]{Corollary}
\newtheorem{lemma}[theorem]{Lemma}
\newtheorem{proposition}[theorem]{Proposition}
\newtheorem{axiom}[theorem]{Axiom}
\newtheorem{conjecture}[theorem]{Conjecture}
\newtheorem{definition}[theorem]{Definition}
\newtheorem{remark}[theorem]{Remark}
\newtheorem{notation}[theorem]{Notation}
\newtheorem{exercise}[theorem]{Exercise}
\newtheorem{example}[theorem]{Example}
\numberwithin{equation}{section}

\section{Introduction}

The Euler characteristic is one of the simplest homotopy invariants of a
smooth, closed manifold. We begin by briefly reviewing some standard
theorems of topology which establish the equivalence of three ways of
computing it. Then we proceed to new analogs for Riemannian foliations,
where more details will be supplied.

For a smooth, closed manifold $M$, we define the Euler characteristic as 
\begin{equation*}
\chi (M)=\sum (-1)^{k}\dim H^{k}(M),
\end{equation*}
where $H^{k}$ is de Rham cohomology. According to classical results of de
Rham the spaces $H^{k}(M)$ are finite-dimensional and homotopy invariant;
thus the definition of $\chi(M) $ makes sense and is homotopy invariant. An
alternate route to the finite dimensionality of $H^{k}$ and the topological
invariance of $\chi $ goes by way of the \v {C}ech-de Rham double complex,
which computes the de Rham cohomology and provides an isomorphism with the
\v {C}ech cohomology of a finite good cover of $M$ (here ``good'' means that
non-empty finite intersections are diffeomorphic to $\Bbb{R}^{n}$). Using
this isomorphism one can define the Euler class $e(E)\nolinebreak \in
\nolinebreak H^{n}(M)$ of an oriented $(n-1)$-sphere bundle $E\rightarrow M$
. The vanishing of the Euler class is a necessary, but not in general
sufficient, condition for the bundle to have a section. In the case that $e$
is the Euler class of the sphere bundle associated with the tangent bundle
of $M$, $\int_{M}e=\chi (M)$. A corollary of this is the Hopf Index Theorem: 
$\chi (M)$ is the sum of the indices of the singular points of any
non-degenerate vector field on $M$. Details can be found in \cite{BT}.

Suppose that a smooth, closed manifold $M$ is endowed with a smooth
foliation $\mathcal{F}$. A form $\omega $ on $M$ is \emph{basic} if for
every vector field $X$ tangent to the leaves, $i(X)\omega =0$ and $
i(X)(d\omega )=0$, where $i(X)$ denotes interior product with $X$. The
exterior derivative of a basic form is again basic, so the basic forms are a
subcomplex $\Omega _{B}^{\ast }\left( M,\mathcal{\ F}\right) $ (or $\Omega
_{B}^{\ast }\left( M\right) $) of the de Rham complex $\Omega ^{\ast }\left(
M\right) $. The cohomology of this subcomplex is the \emph{basic cohomology} 
$H_{B}^{\ast }\left( M,\mathcal{F}\right) .$

Suppose $\mathcal{F}$ has codimension $q$. The \emph{basic Euler
characteristic} is defined as 
\begin{equation*}
\chi _{B}\left( M,\mathcal{F}\right) =\sum_{k=0}^{q}(-1)^{k}\dim
H_{B}^{k}\left( M,\mathcal{F}\right) ,
\end{equation*}
provided that all of the basic cohomology groups are finite-dimensional.
Although $H_{B}^{0}\left( M,\mathcal{F}\right) $ and $H_{B}^{1}\left( M, 
\mathcal{F}\right) $ are always finite-dimensional, there are foliations for
which higher basic cohomology groups can be infinite-dimensional. For
example, in \cite{Gh}, the author gives an example of a flow on a 3-manifold
for which $H_{B}^{2}\left( M,\mathcal{F}\right) $ is infinite-dimensional.
Therefore, we must restrict our investigation to a class of foliations for
which the basic cohomology is finite-dimensional. A large class of
foliations with this property are \emph{Riemannian foliations}; a foliation
is Riemannian if its normal bundle admits a holonomy-invariant Riemannian
metric. There are various proofs that the basic cohomology of a Riemannian
foliation on a closed manifold is finite-dimensional; see for example \cite
{EKHS} for the original proof using spectral sequence techniques or \cite{KT}
and \cite{PR} for proofs using a basic version of the Hodge theorem.

In Section 2 we develop a basic version of \v {C}ech-de Rham cohomology
along the lines of \cite{BT}. The assumption that the foliation is
Riemannian is needed to obtain basic partitions of unity and the basic
Mayer-Vietoris sequence. For the basic Poincar\'{e} lemma we further assume
that all of the leaves are closed. Examples are given to show that these
conditions are necessary. The basic \v {C}ech-de Rham theorem establishes
the equivalence of $H_{B}^{*}$ and $\check{H}_{B}^{*}$, if all the leaves
are closed. Examples show that the basic \v {C}ech cohomology and the basic
de Rham cohomology are not necessarily isomorphic for Riemannian foliations
in general.

A primary goal of this paper is to establish a foliation version of the Hopf
Index Theorem. However, the standard proofs of this theorem do not carry
over, even for Riemannian foliations. The problem is there are many
Riemannian foliations that have a nonzero basic Euler characteristic, yet
have trivial top-dimensional basic cohomology (for example, a non-taut
Riemannian foliation (\cite{To}) can possess this property). Thus it is
impossible for these foliations to have any sort of basic Euler class that
can be integrated to obtain the basic Euler characteristic (we do mention
that for \emph{taut} Riemannian foliations, one can define a nontrivial
basic Euler class \cite{Z}).

To establish a Hopf Index Theorem for Riemannian foliations, another
approach is required. The approach we use, which is in Section 3, involves a
modification of the Witten deformation of the de Rham complex to the
foliation case. Let $V$ be a basic vector field on $(M,\mathcal{F})$; i.e.,
a vector field on $M$ whose flow maps leaves to leaves (see \cite{Mo} or 
\cite{To}, for example). In analogy with what is needed in the classical
Hopf Index Theorem, we require that $V$ satisfy a transverse nondegeneracy
condition, which we call $\mathcal{F}$-nondegeneracy. We say that a leaf
closure is \emph{critical} for $V$ if $V$ is tangent to $\mathcal{F}$ at
that leaf closure. When $V$ is $\mathcal{F}$-nondegenerate, the critical
leaf closures are necessarily isolated. For each critical leaf closure $L$,
we define the \emph{index} ind$_{L}(V)$ of $V$ at $L$; just as in the
classical case, this index is always $\pm 1$.

Next, let $d_B$ be the restriction of the exterior derivative to basic
forms, let $\delta_B$ denote the formal adjoint of $d_B$, and for each real
number $s$, define the basic Witten differential

\begin{equation*}
D_{B,s} = d_{B}+s\,i\left( V\right) + \delta _{B}+sV^{\flat }\wedge :\Omega
_{B}^{*}\left( M\right) \rightarrow \Omega _{B}^{*}\left( M\right).
\end{equation*}

We show that the index of $D_{B,s}$ is independent of $s$, and examine the
behavior of this operator as $s$ goes to infinity. In the limit, the formula
for the index of $D_{B,s}$ concentrates at the critical leaf closures. We
next establish the necessary analytic properties of the basic Witten
deformation. This leads to Theorem \ref{BasicHopf}: Let $\left( M,\mathcal{F}
\right) $ be a Riemannian foliation, and let $V$ be a basic vector field
that is $\mathcal{F}$-nondegenerate. Given a critical leaf closure $L$, let $
\mathcal{\ O}_{L}=$ $\mathcal{O}_{L}\left( V\right)$ denote the orientation
line bundle of $V$ at $L$ (Definition \ref{orientation}). Then 
\begin{equation*}
\chi _{B}\left( M,\mathcal{F}\right) =\sum_{L\,\mathrm{critical} }
\mathrm{%
ind}_{L}\left( V\right) \chi _{B}\left( L,\mathcal{%
F},\mathcal{O}_{L}\right) .
\end{equation*}
where $\chi _{B}\left( L,\mathcal{F},\mathcal{O}_{L}\right)$ is the
alternating sum of the dimensions of the certain cohomology groups $
H_B^*\left( L,\mathcal{F},\mathcal{O}_{L}\right)$ associated to the
foliation $\mathcal{F}$ restricted to $L$(Definition \ref{twisted}). We
remark that in many simple cases, $\chi
_{B}\left( L,\mathcal{F}, \mathcal{O}_{L}\right) = 1$, whence our formula
takes on the precise form of the classical Hopf Index Theorem.

One important implication (Corollary~\ref{obstruction}) of this result is
that if $(M,\mathcal{F})$ admits a basic vector field that is never tangent
to the leaves, then $\chi _{B}\left( M,\mathcal{F}\right) = 0$. For example,
this implies that if $\mathcal{F}$ has codimension 1 or 2 and has such a
basic vector field, then $M$ has infinite fundamental group (Corollary~\ref
{fundgroup}).

\section{Basic \v{C}ech and de Rham Cohomologies}

Let $M$ be a compact manifold, and let $\mathcal{F}$ be a Riemannian
foliation on $M$. Let $J$ be a finite ordered set, and let $\mathcal{U}
=\left\{ U_{\alpha }\right\} _{\alpha \in J}$ be a finite basic open cover
of $M$; \textit{i.e.}, a finite open cover in which each set $U_{\alpha }$
is a union of leaves. Such an open cover always exists in this case, because
tubular neighborhoods of leaf closures are unions of leaves. For $\alpha
_{0},\alpha _{1},\cdots ,\alpha _{n}$ an increasing sequence of indices
define $U_{\alpha _{0}\alpha _{1}\cdots \alpha
_{n}}=\bigcap_{i=0}^{n}U_{\alpha _{i}}$, and let $\Omega _{B}^{p}(U_{\alpha
_{0}\alpha _{1}\cdots \alpha _{n}})$ be the collection of basic $p$-forms on 
$U_{\alpha _{0}\alpha _{1}\cdots \alpha _{n}}$. Define 
\begin{equation*}
\delta :\prod \Omega _{B}^{*}(U_{\alpha _{0}\alpha _{1}\cdots \alpha
_{k}})\longrightarrow \prod \Omega _{B}^{*}(U_{\alpha _{0}\alpha _{1}\cdots
\alpha _{k+1}})
\end{equation*}
by the formula 
\begin{equation*}
(\delta \omega )_{\alpha _{0}\alpha _{1}\cdots \alpha
_{k+1}}=\sum_{i=0}^{k+1}(-1)^{i}\omega _{\alpha _{0}\cdots \widehat{\alpha }
_{i}\cdots \alpha _{k+1}}\text{,}
\end{equation*}
where each form is restricted to the appropriate subset. It is convenient to
extend the index set of the product to include \emph{all} sequences of the
appropriate length from $J$, regardless of order or repetition, and to adopt
the convention that interchanging indices of a component introduces a minus
sign. So, for example, $\omega _{\alpha _{2}\alpha _{1}}=-\omega _{\alpha
_{1}\alpha _{2}}$ for all $\alpha _{1}$and $\alpha _{2}$ in $J$. It is
straightforward to check that the definition of $\delta $ respects the sign
convention and $\delta ^{2}=0$.

\begin{definition}
A \textrm{basic partition of unity subordinate to a basic open cover} $
\left\{ U_{\alpha }\right\} $ is a partition of unity $\left\{ \rho _{\alpha
}\right\} $ consisting of basic functions $\rho _{\alpha }$; that is,
functions that are constant on leaves.
\end{definition}

\begin{lemma}
Every basic open cover of a Riemannian foliation admits a basic partition of
unity.
\end{lemma}

\begin{proof}
Endow the manifold with a bundle-like metric, and choose \emph{any}
partition of unity subbordinate to the basic cover. Orthogonally project the
functions in the partition of unity to the space of basic functions; the
smoothness and other desired properties of the resulting functions are
guaranteed by the results of \cite{PR}.
\end{proof}

\begin{theorem}
(Basic Mayer-Vietoris Sequence) \label{BMV}The sequence 
\begin{equation*}
0\longrightarrow \Omega _{B}^{\ast }(M)\overset{r}{\longrightarrow }\prod
\Omega _{B}^{\ast }(U_{\alpha _{0}})\overset{\delta }{\longrightarrow }\prod
\Omega _{B}^{\ast }(U_{\alpha _{0}\alpha _{1}})\overset{\delta }{%
\longrightarrow }\cdots
\end{equation*}
is exact, where $r$ denotes the restriction map.
\end{theorem}

\begin{proof}
The proof of this theorem is the same as the proof of Theorem II.8.5 in \cite
{BT}, with forms replaced by basic forms and partitions of unity replaced by
basic partitions of unity.
\end{proof}

\begin{remark}
If the foliation is not Riemannian, it may not admit basic partitions of
unity, and in such cases, the Mayer-Vietoris sequence may fail to be exact.
To describe such an example, we begin with some notation. For each real
number $t$ and open interval $\left( a,b\right) $, let $L_{\left( a,b\right)
}^{t}$ be the curve $\left\{ \left( x,\frac{1}{x-a}+\frac{1}{b-x}+t\right)
\,\left| \,a<x<b\right. \right\} $ in $\Bbb{R}^{2}$. Foliate $\lbrack
0,1)\times \Bbb{R}$ by the lines $x=$ $\frac{1}{3}$ and $x=$ $\frac{2}{3},$
along with Reeb components $L_{\left( 0,\frac{1}{3}\right) }^{t},L_{\left( 
\frac{1}{3},\frac{2}{3}\right) }^{t}$, and $L_{\left( \frac{2}{3},1\right)
}^{t}$, where $t$ ranges over all real numbers. Let $\Bbb{T}^{2}\subset \Bbb{%
C}^{2}$ be the 2-torus, and consider the open map $p:\lbrack 0,1)\times \Bbb{%
R\longrightarrow T}^{2}$ given by $p(x,y)=(e^{2\pi ix},e^{2\pi iy})$. Then
the image under $p$ of our foliation on $\lbrack 0,1)\times \Bbb{R}$
determines a foliation on $\Bbb{T}^{2}$. Let $U_{0}=\left( 0,\frac{2}{3}
\right) \times \Bbb{R}$, $U_{1}=\left( \frac{1}{3},1\right) \times \Bbb{R}$
,and $U_{2}=\left( \left[ 0,\frac{1}{3}\right) \cup \left( \frac{2}{3}
,1\right) \right) \times \Bbb{R}$. There does not exist a basic partition of
unity subordinate to the basic open cover $\left\{ p\left( U_{0}\right)
,p\left( U_{1}\right) ,p\left( U_{2}\right) \right\} $ of $\Bbb{T}^{2}$,
because the only basic functions are the constant ones. To see that the
Mayer-Vietoris sequence is not exact in this example, first note that the
intersection of all three sets in our basic open cover is empty. Therefore,
for basic functions, we have the following piece of the Mayer-Vietoris
sequence: 
\begin{multline*}
C_{B}^{\infty }\left( p\left( U_{0}\right) \right) \oplus C_{B}^{\infty
}\left( p\left( U_{1}\right) \right) \oplus C_{B}^{\infty }\left( p\left(
U_{2}\right) \right) \overset{\delta }{\longrightarrow } \\
C_{B}^{\infty }\left( p\left( U_{0}\cap U_{1}\right) \right) \oplus
C_{B}^{\infty }\left( p\left( U_{1}\cap U_{2}\right) \right) \oplus
C_{B}^{\infty }\left( p\left( U_{0}\cap U_{2}\right) \right) \longrightarrow
0.
\end{multline*}

The map $\delta $ here is given by the formula $\delta
(c_{0},c_{1},c_{2})=(c_{1}-c_{0},c_{2}-c_{1},c_{2}-c_{0})$, which is not
surjective since there is a nontrivial relation among the three components
of the image of $\delta $.
\end{remark}

\begin{remark}
Even if a basic open cover of a (necessarily) non-Riemannian foliation does
not admit basic partitions of unity, it is possible for the Mayer-Vietoris
sequence to be exact. One such example is as follows; we maintain the
notation of the previous remark. Foliate $\lbrack 0,1)\times \Bbb{R}$ by
leaves of two sorts: straight lines $L_{a}=\left\{ \left( a,y\right)
\,|\,y\in \Bbb{R}\right\} $ for $0\leq a\leq \frac{1}{2}$, and a Reeb
component $L_{\left( \frac{1}{2},1\right) }^{t},t\in \Bbb{R}$. As in the
previous remark, the image under $p$ of this foliation of $\lbrack
0,1)\times \Bbb{R}$ determines a foliation on $\Bbb{T}^{2}.$ Let $
U_{0}=(0,1)\times \Bbb{R}$ and $U_{1}=\left( \left[ 0,\frac{1}{4}\right)
\cup \left( \frac{1}{2},1\right) \right) \times \Bbb{R}$. Then there does
not exist a basic partition of unity subordinate to the basic open cover $
\left\{ p\left( U_{0}\right) ,p\left( U_{1}\right) \right\} $ of $\Bbb{T}%
^{2} $, because every basic function on the Reeb component must be constant.
In spite of this, it is straightforward to check that the Mayer-Vietoris
sequence associated to this basic open cover is exact.
\end{remark}

There are various cohomology theories that one can associate to the
foliation $\mathcal{F}\frak{.}$ Note that the exterior derivative $d\omega $
of a basic form $\omega$ is basic, and therefore the collection of basic
forms of $\mathcal{F}$ is a subcomplex of the de Rham complex. The
cohomology of this subcomplex is denoted $H_{B}^{*}(M,\mathcal{F}\frak{)}$,
and is called the \textit{basic de Rham cohomology of }$(M,\mathcal{F}\frak{)%
}$.

Another cohomology theory can be defined for any basic open cover $\mathcal{U%
}$ of $(M,\mathcal{F}\frak{)}$. For each $p\geq 0$, define 
\begin{equation*}
C_{B}^{p}(\mathcal{U},\Bbb{R})=\text{ker}\left( d:\prod \Omega
_{B}^{0}(U_{\alpha _{0}\cdots \alpha _{p}})\longrightarrow \prod \Omega
_{B}^{1}(U_{\alpha _{0}\cdots \alpha _{p}})\right) .
\end{equation*}
It is easy to check that $\left( C_{B}^{*}(\mathcal{U},\Bbb{R}),\delta
\right) $ is a cochain complex; we declare its cohomology $\check{H}
_{B}^{p}( \mathcal{U},\Bbb{R})$ to be the \textit{basic \v {C}ech cohomology}
of the basic open cover $\mathcal{U}$.

Our immediate goal is to show that for ``nice'' covers $\mathcal{U}$, basic
de Rham cohomology and basic \v {C}ech cohomology are isomorphic. To do
this, we need a foliation version of the Poincar\'{e} lemma.

\begin{proposition}
(Basic Poincar\'{e} Lemma) Let $\mathcal{F}$ be a Riemannian foliation of $M$
, suppose that all of the leaves of $\mathcal{F}$ are closed, and equip $M$
with a bundle-like metric. Let $L$ be a leaf; for $\varepsilon >0$, let $U$
be the tubular neighborhood of $L$ consisting of points that are a distance
less than $\varepsilon $ from $L.$ Then, for $\varepsilon $ sufficiently
small and for $k>0$, every closed basic $k$-form on $U$ is exact.
\end{proposition}

\begin{proof}
First observe that since we have given $M$ a bundle-like metric, the tubular
neighborhood $U$ is a union of leaves, and so the restriction $\mathcal{F}
\left| _{U}\right. $ of $\mathcal{F}$ to $U$ makes sense. Choose $
\varepsilon $ small enough so that $U$ misses the cut locus of $L$; since $M$
is compact and $L$ is closed, this can always be done. Fix $x\in L$, and let 
$D$ be the exponential image of the ball of radius $\varepsilon $ in the
normal space $N_{x}L$. Then 
\begin{equation*}
\Omega _{B}^{k}(U)\cong \left\{ \eta \in \Omega ^{k}(D)\,|\,\,\eta \,\text{
is holonomy invariant}\right\} \text{.}
\end{equation*}
The holonomy of $\mathcal{F}$ acts by a finite subgroup $\Gamma $ of the
orthogonal group (\cite{Mo}), so we have 
\begin{equation*}
\Omega _{B}^{k}(U)\cong \left\{ \eta \in \Omega ^{k}(D)\left| g^{*}\eta
=\eta \right. \text{for all }g\text{ in }\Gamma \right\} \text{,}
\end{equation*}
and this isomorphism commutes with the exterior derivative. Suppose $\omega
\nolinebreak \in \nolinebreak \Omega _{B}^{k}(U)$ is closed, and let $\eta
\nolinebreak \in \nolinebreak \Omega ^{k}(D)$ be the closed form associated
to $\omega $ via the isomorphism above. Since $D$ is diffeomorphic to
Euclidean space, there exists by the ordinary Poincar\'{e} lemma a form $\mu
\nolinebreak \in \nolinebreak \Omega ^{k-1}(D)$ such that $d\mu =\eta $.
Now, $\mu $ may not be $\Gamma $-invariant, but the averaged form $\zeta = 
\frac{1}{|\Gamma |}\sum_{g\in \Gamma }g^{*}\mu $ is, and $d\zeta $ also
equals $\eta $. Therefore $\eta $, and hence $\omega $, is exact.
\end{proof}

\begin{definition}
A \textrm{basic good cover of }$(M,\mathcal{F}\frak{)}$ is a basic open
cover $\mathcal{U}$ of $(M,\mathcal{F}\frak{)}$ with the feature that the
basic cohomology of each intersection is trivial.
\end{definition}

When all of the leaves of $\mathcal{F}$ are closed, the basic Poincar\'{e}
lemma implies that we can obtain a basic good cover by covering each leaf by
a sufficiently small tubular neighborhood. In fact, since $M$ is compact, we
can choose a \textit{finite }basic good cover of $M.$

\begin{theorem}
(Basic \v{C}ech-de Rham Theorem) Suppose $\mathcal{F}$ is Riemannian and
consists entirely of closed leaves, and suppose that $\allowbreak \mathcal{U}
$ is a basic good cover of $M$. Then 
\begin{equation*}
H_{B}^{\ast }\left( M,\mathcal{F}\right) \frak{\cong }\check{H}_{B}^{\ast
}\left( \mathcal{U},\Bbb{R}\right) .
\end{equation*}
\end{theorem}

\begin{proof}
For each $p,q\geq 0$, let $C^{p}\left( \mathcal{U},\Omega _{B}^{q}\right) $
be the collection of basic $q$-forms restricted to $\left( p+1\right) $-fold
intersections, and define a double complex 
\begin{equation*}
C^{*}\left( \mathcal{U},\Omega _{B}^{*}\right) =\bigoplus_{p,q\geq
0}C^{p}\left( \mathcal{U},\Omega _{B}^{q}\right) .
\end{equation*}
The horizontal differential is the map $\delta $ defined above, and the
vertical differential is $(-1)^{p}d$, where $d$ is the exterior derivative.
We compute the cohomology {$H^{*}\left( \mathcal{U},\Omega _{B}^{*}\right) $}
of this double complex in two ways.

First, augment $C^{*}\left( \mathcal{U},\Omega _{B}^{*}\right) $ by the
column $\bigoplus_{q\geq 0}\Omega _{B}^{q}(M)$, and then map into \linebreak$
\bigoplus_{q\geq 0}C^{1}\left( \mathcal{U},\Omega _{B}^{q}\right) $ by the
restriction map $r$. By Theorem~\ref{BMV}, the rows of the augmented double
complex are exact. By \cite[p. 97]{BT}, the cohomology $H^{*}\left( \mathcal{%
U},\Omega _{B}^{*}\right) $ of the original double complex is isomorphic to
the cohomology of the initial column, which is precisely the basic de Rham
cohomology $H_{B}^{*}(M,\mathcal{F}\frak{)}$ of $(M,\mathcal{F}\frak{)}$.

Second, augment the original double complex with the row $\bigoplus_{p\geq
0}C_{B}^{p}\left( \mathcal{U},\Bbb{R}\right) $. Then the columns of the new
double complex are exact, and thus $\check{H}_{B}^{*}\left( \mathcal{U},\Bbb{%
R}\right) $ is also isomorphic to $H^{*}\left( \mathcal{U},\Omega
_{B}^{*}\right) $.
\end{proof}

\begin{corollary}
The group $\check{H}_{B}^{p}\left( \mathcal{U},\Bbb{R}\right) $ is
independent of the choice of basic good cover $\mathcal{U}$.
\end{corollary}

\begin{corollary}
If $\mathcal{F}$ is a Riemannian foliation and consists entirely of closed
leaves, then $H_{B}^{\ast }(M,\mathcal{F}\frak{)}$ is finite-dimensional.
\end{corollary}

In fact, the group $H_{B}^{*}(M,\mathcal{F}\frak{)}$ is finite-dimensional
for any Riemannian foliation $\mathcal{F}$, as mentioned in the introduction.

\begin{remark}
The basic \v{C}ech-de Rham theorem is not necessarily true if the Riemannian
foliation $\mathcal{F}$ has some non-closed leaves. Consider the foliation
of the 2-torus $\Bbb{T}^{2}$ foliated by translates of a line of irrational
slope. In this case, the only nonempty basic open set is $\Bbb{T}^{2}$
itself, so there is only one basic open cover $\mathcal{U}$ to choose.
Clearly $\check{H}_{B}^{1}\left( \mathcal{U},\Bbb{R}\right) \cong $ $0,$
whereas it is straightforward to check that $H_{B}^{1}(M,\mathcal{F}\frak{)}$
$\cong \Bbb{R}$. Note also that since 
$H_{B}^{1}(M,\mathcal{F})$ is nontrivial, 
the basic Poincar\'{e} lemma fails in this example.
\end{remark}

\section{Witten deformation of the basic de Rham complex}

For the remainder of this paper, $\mathcal{F}$ is a Riemannian foliation on
a smooth compact manifold $M$, and $M$ is equipped with a bundle-like
metric; $\omega _{\mathrm{vol}}$ will denote the volume form associated to
the metric.

Let $V$ be a smooth vector field on $\left( M,\mathcal{F}\right) $. We say
that $V$ is a\emph{\ basic vector field} if for every vector field $X$ in $T 
\mathcal{F}$, $\left[ V,X\right] \ $is in $T\mathcal{F}$. If in addition $
V\left( x\right) $ is in $N_{x}\mathcal{F}$ (the normal space to the leaf at 
$x$) for every $x\in M,$ we will say that $V$ is a \emph{basic normal vector
field}. Basic normal vector fields always exist; the projection of any basic
vector field onto $N\mathcal{F}$ is such a vector field. Associated to a
basic vector field $V$ is a one-parameter family of diffeomorphisms of $M$
that preserves the foliation $\mathcal{F}$.

Let $V$ be a basic vector field. Let $L$ be a leaf closure with the property
that $V$ is tangent to every leaf in $L$; such a leaf closure will be called
a \emph{critical leaf closure} for $V$. We note in passing that if a basic
vector field is tangent to the foliation at any point, it is in fact tangent
to every leaf in the leaf closure. Define the \emph{linear part} of $V$ at $
x\in L$ to be the linear map $V_{L}:N_{x}L\rightarrow $ $N_{x}L$ defined by $
X\mapsto \pi \left[ V,\overline{X}\right] _{x}$, where $\overline{X}$ is any
vector field that restricts to $X$ at $x$, and $\pi :T_{x}M\rightarrow $ $
N_{x}L$ is the projection map. The basic vector field $V$ will be called $ 
\mathcal{F}$\emph{-nondegenerate} if the linear part of $V$ is an
isomorphism at each point of every critical leaf closure. Every Riemannian
foliation $\mathcal{F}$ admits $\mathcal{F}$-nondegenerate\emph{\ }basic
vector fields; simple examples are gradients of basic Morse functions (see 
\cite{A}). A critical leaf closure for a $\mathcal{F}$-nondegenerate basic
vector field $V$ is necessarily isolated. We say that a nondegenerate vector
field $V$ has \emph{index }$1$ (respectively, \emph{\ index }$-1$) at a
critical leaf closure $L$ if the determinant of the linear transformation $
V_{L}$ is everywhere positive (respectively, negative) on $L$. Clearly, if $
V $ is $\mathcal{F}$-nondegenerate, then at each critical leaf closure, $V$
must have either index $1$ or $-1$.

Given any point $x_{0}$ of a critical leaf closure $L$, choose orthonormal
coordinates $\overline{y}=\left( y_{1},...,y_{\bar{q}}\right) $ for the
normal space $N_{x}\overline{\mathcal{F}}=N_{x}L$, and extend these
coordinates to orthonormal coordinates $y=\left( y_{1},...,y_{\bar{q}},y_{%
\bar{q}+1},...,y_{q}\right) $ for $N_{x}\mathcal{F}$. Let $x=\left(
x_{1},...,x_{p}\right) $ be geodesic normal coordinates for the leaf near $
x_{0}$. The coordinates $\left( x,y\right) $ parametrize a tubular
neighborhood of $L$ near $x_{0}$ via the normal exponential map. It is
elementary to check that we may write $V$ near $x_{0}$ as an orthogonal sum $
V=V_{1}+V_{2}+V_{3},$ where $V_{1}$ is tangent to $\mathcal{F}$, $
V_{2}\left( y\right) =\sum_{i=1}^{\bar{q}}\alpha ^{i}\left( y\right) \frac{%
\partial }{\partial y_{i}}$, and $V_{3}\left( y\right) =\sum_{j=_{\bar{q}%
+1}}^{q}\beta ^{j}\left( y\right) \frac{\partial }{\partial y_{j}}$. The
linear transformation $V_{L}$ is given by multiplication by the $\bar{q}
\times \bar{q}$ matrix $\left( \frac{\partial \alpha ^{i}}{\partial y_{m}}
\left( y\right) \right) _{1\leq i,m\leq \bar{q}}$ on $N_{x}L$, so that the
index of $V$ is simply the sign of the determinant of that matrix. We remark
that if $V$ is the gradient of a basic function $f$, then $V_{L}$ is the
Hessian of $f$ restricted to the normal ball, and $V_{1}=V_{3}=0.$

Observe that a basic normal vector field $V$ (hence $V_{1}=0$) is determined
on a tubular neighborhood $U$ of the leaf closure $L$ by its values on a
transverse $\bar{q}$-dimensional ball. Let $B\left( \delta \right) $ be the
image of a ball of radius $\delta $ under the normal exponential map $\exp
_{x}^{\bot }:N_{x}L\rightarrow U$ at $x\in L$, so that $y=\left(
y_{1},...,y_{\bar{q}}\right) $ are exponential coordinates for this ball.
The holonomy group at $x$ is represented on $B\left( \delta \right) $ by a group $
\Gamma $ of orthogonal transformations (see \cite{Mo}), so that we may
identify basic normal vector fields on $U$ with normal vector fields
restricted to $B\left( \delta \right) $ that are equivariant with respect to
the action of $\Gamma $. That is, for any $g\in \Gamma $, every basic normal
vector field restricted to $B\left( \delta \right) $ has the form 
\begin{equation*}
V\left( y\right) =\sum_{i=1}^{q}v^{i}\left( \overline{y}\right) \frac{%
\partial }{\partial y_{i}}=\left( 
\begin{array}{c}
v^{1}\left( \overline{y}\right) \\ 
\vdots \\ 
v^{q}\left( \overline{y}\right)
\end{array}
\right)
\end{equation*}
and satisfies 
\begin{equation*}
g\left( 
\begin{array}{c}
v^{1}\left( \overline{y}\right) \\ 
\vdots \\ 
v^{q}\left( \overline{y}\right)
\end{array}
\right) =\left( 
\begin{array}{c}
v^{1}\left( g\overline{y}\right) \\ 
\vdots \\ 
v^{q}\left( g\overline{y}\right)
\end{array}
\right) .
\end{equation*}
Note that the notation $g\overline{y}$ denotes the restriction of the action
of $g$ on $N_{x}\mathcal{F}$ to the subspace $N_{x}L$. This makes sense,
because elements of the holonomy group map vector fields orthogonal to $L$
to other vector fields orthogonal to $L$. As a consequence, the action of $
\Gamma $ commutes with the action of $V_{L}$ on $N_{x}L$.

\begin{definition}
\label{orientation}Let $V$ be a basic vector field that has a nondegenerate
zero on a critical leaf closure $L$. At each point $x\in L$, let $
V_{L,x}:N_{x}L\rightarrow N_{x}L$ be the linear part of $V$ at $x$. Since $
V_{L,x}$ is invertible, it can be written as the (uniquely determined) polar
decomposition $V_{L,x}=P_{x}\Theta _{x}$, where $P_{x}=\sqrt{V_{L,x}^{\ast
}V_{L,x}}$ is positive and symmetric, and $\Theta _{x}=P_{x}^{-1}V_{L,x}$ is
an isometry. The \textrm{orientation line bundle }$\mathcal{O}_{L}=\mathcal{O%
}_{L}\left( V\right) $ \textrm{of }$V$\textrm{\ at }$L$ is the orientation
bundle of the subbundle of $NL$ spanned by the eigenvectors of $\Theta _{x}$
corresponding the eigenvalue $-1$.
\end{definition}

\begin{remark}
It can be shown that $\mathcal{O}_{L}$ is a smooth line bundle. Observe that
in the context of this paper, the singularity of the subbundle is not an
issue, since the eigenvalues of $\Theta _{x}$ are constant as $x$ moves
along the leaves. We also observe that if $V=\nabla f$ for a Bott-Morse
function $f$, then this orientation line bundle corresponds to the
orientation line bundle of the subbundle of negative directions of $f$.
\end{remark}

We now prove a result that puts a basic vector field into standard form.

\begin{lemma}
\label{deform}Let $V$ be a basic vector field that is $\mathcal{F}$
-nondegenerate, and let $L$ be a critical leaf closure. Let $T_{\delta
}\left( L\right) $ denote the tubular neighborhood of radius $\delta $
around $L$, and assume that $\delta $ is chosen so that $T_{\delta }\left(
L\right) $ does not contain any other critical leaf closures and misses the
cut locus of $L$. Then there is a basic normal vector field $\widetilde{V}$
and a $\widetilde{\delta }$ with $0<\widetilde{\delta }<\delta $ such that:

\begin{enumerate}
\item  $\widetilde{V}=V$ outside $T_{\delta }\left( L\right) $;

\item  $L$ is the only critical leaf closure of $\widetilde{V}$ in $
T_{\delta }\left( L\right) $;

\item  The index of $\widetilde{V}$ at $L$ is the same as the index of $V$
at $L$;

\item  $\widetilde{V}=\nabla f$ on $T_{\widetilde{\delta }}\left( L\right) $
, where $f$ is a basic function whose Morse index (restricted to a normal
ball at a point of $L$) is even if the index of $V$ is $1$ at $L$ and is odd
if the index of $V$ is $-1$ at $L$;

\item  The orientation line bundle of $\widetilde{V}$ at $L$ is identical to
the orientation line bundle of $V$ at $L$.
\end{enumerate}
\end{lemma}

\begin{proof}
Choose coordinates $\left( x,y\right) =\left( x_{1},...,x_{p},y_{1},...,y_{%
\bar{q}},y_{\bar{q}+1},...,y_{q}\right) $ near a point $x\in L$ as described
in the paragraphs above. As before, we write $V$ in the form $
V=V_{1}+V_{2}+V_{3},$ where $V_{1}$ is tangent to $\mathcal{F}$, $
V_{2}\left( y\right) =\sum_{i=1}^{\bar{q}}\alpha ^{i}\left( y\right) \frac{%
\partial }{\partial y_{i}}$, and $V_{3}\left( y\right) =\sum_{j=_{\bar{q}%
}+1}^{q}\beta ^{j}\left( y\right) \frac{\partial }{\partial y_{j}}$. Given
any $\delta ^{\prime }$ such that $0<\delta ^{\prime }<\delta $, we may
multiply the tangent component $V_{1}$ by a radial basic function that is
zero in $T_{\delta ^{\prime }}\left( L\right) $ and is $1$ outside $
T_{\delta }\left( L\right) $. In doing so, we preserve the index and
orientation line bundle and yet restrict to the case where $V$ is a basic
normal vector field. We now assume that $V$ has already been modified in
this way. Next, the basic normal vector field $V$ is determined on $
T_{\delta ^{\prime }}\left( L\right) $ by its restriction to a normal ball $
B\left( \delta ^{\prime }\right) $ with the $\overline{y}=\left(
y_{1},...,y_{\bar{q}}\right) $ coordinates, and we may write $V\left( 
\overline{y}\right) =\sum_{i=1}^{\bar{q}}v^{i}\left( \overline{y}\right) 
\frac{\partial }{\partial y_{i}}+\sum_{j=\bar{q}+1}^{q}v^{j}\left( \overline{%
y}\right) \frac{\partial }{\partial y_{j}}$ on $B\left( \delta ^{\prime
}\right) $. Again, we may multiply the component $\sum_{j=\bar{q}%
+1}^{q}v^{j}\left( \overline{y}\right) \frac{\partial }{\partial y_{j}}$ by
a similar radial basic function if necessary so that this component vanishes
on a given $B\left( \delta ^{\prime \prime }\right) $ such that $0<\delta
^{\prime \prime }<\delta ^{\prime }$, without changing the relevant
properties of $V$. For the remainder of this proof, assume that we have
already modified $V$ so that it is in the form $V=\sum_{i=1}^{\bar{q}%
}v^{i}\left( \overline{y}\right) \frac{\partial }{\partial y_{i}}$
restricted to $B\left( \delta \right) $, choosing a slightly smaller $\delta 
$ if necessary. It follows that $V\left( \overline{y}\right) =V_{L}\left( 
\overline{y}\right) +O\left( \left\| \overline{y}\right\| ^{2}\right) $ on $
B\left( \delta \right) $. Given $x\in L$, we write $V_{L}:N_{x}L\rightarrow
N_{x}L$ in terms of its polar decomposition $V_{L}=P\Theta $. Let $\Gamma $
be the closed subgroup of isometries on $N_{x}L$ induced from the
representation of the holonomy group at $x$ on $N_{x}L$. For every $g\in
\Gamma $, $g\,V_{L}=V_{L}\,g$, whence $V_{L}^{\ast
}\,g^{-1}=g^{-1}V_{L}^{\ast }$. This in turn implies that every $g\in \Gamma 
$ commutes with $P$ and $\Theta $. Let $P_{t}$ be defined by $
P_{t}v_{i}=\lambda _{i}\left( t\right) v_{i}$ on the $\lambda _{i}$
-eigenspace of $P$, where each $\lambda _{i}\left( t\right) $ is any smooth
positive function such that $\lambda _{i}\left( 0\right) =\lambda _{i}$ and $
\lambda _{i}\left( 1\right) =1$. Since every $g\in \Gamma $ and $P$ have
simultaneous eigenspaces, every $g\in \Gamma $ also commutes with each $
P_{t} $. Thus, the smooth, one-parameter family of transformations $\left\{
T_{t}=P_{t}\Theta \right\} $ is a deformation of $V_{L}$ ($t=0$) to an
orthogonal transformation $\Theta $ ($t=1$) that has constant index (that
is, the sign of the determinant of the linear transformation does not
change). Next, since $\Theta $ is orthogonal, there is a complex orthogonal
basis $\left\{ w_{k}\right\} $ of $N_{x}L\otimes \Bbb{C}$ consisting of
eigenvectors such that $\Theta w_{k}=e^{i\theta _{k}}w_{k},$ where $0\leq
\theta _{k}<2\pi $. If $\theta _{k}$ is $0$, then $\Theta $ acts by the
identity on $\mathrm{span}\left\{ \mathrm{Re}w_{k},\mathrm{Im}w_{k}\right\} $
. If $\theta _{k}$ is $\pi $, then $\Theta $ multiplies each vector in $ 
\mathrm{span}\left\{ \mathrm{Re}w_{k},\mathrm{Im}w_{k}\right\} $ by $-1$. If 
$\theta _{k}\neq 0$ and $\theta _{k}\neq \pi $, then $\Theta $ acts by a
rotation of $\theta _{k}$ on $\mathrm{span}\left\{ \mathrm{Re}w_{k},\mathrm{%
Im}w_{k}\right\} $, which in this case is necessarily $2$-dimensional. We
let $\Theta _{0}=\Theta $ and define the transformation $\Theta _{t}$ for $
0<t\leq 1$ by 
\begin{equation*}
\Theta _{t}\left( w_{k}\right) =\left\{ 
\begin{array}{ll}
w_{k} & \text{if }\theta _{k}=0 \\ 
-w_{k} & \text{if }\theta _{k}=\pi \\ 
e^{i\left( 1-t\right) \theta _{k}}w_{k} & \text{otherwise.}
\end{array}
\right.
\end{equation*}
The smooth, one-parameter family of transformations $\left\{ \Theta
_{t}\right\} $ is a deformation of the orthogonal transformation $\Theta $ ($
t=0$) to a transformation $\Theta _{1}$ ($t=1$) that has constant index.
Observe that since each $g\in \Gamma $ commutes with $\Theta _{0},$ the
obvious action of $g$ on $N_{x}L\otimes \Bbb{C}$ satisfies $g\left(
w_{k}\right) =e^{i\alpha _{k}}w_{k}$ for some $\alpha _{k}$; it follows that 
$g$ commutes with each $\Theta _{t}$. The final transformation $\Theta _{1}$
may be described in a real orthogonal basis as a diagonal matrix whose
diagonal consists of $1$'s and $\left( -1\right) $'s; this transformation is
the linear part of a vector field of the form $\sum_{i=1}^{\bar{q}}\pm y^{i} 
\frac{\partial }{\partial y_{i}}$, where the $y_{i}$ are geodesic normal
coordinates on $\exp ^{\perp }\left( N_{x}L\right) $ corresponding to that
particular basis. We also observe that $f\left( \allowbreak \overline{y}
\right) :=\frac{1}{2}\sum_{i=1}^{\bar{q}}\pm \left( y^{i}\right) ^{2}$ is
then a basic function that is well-defined on a small neighborhood of $L$,
and the linear part of $\nabla f\left( \allowbreak \overline{y}\right) $ at $
x\in L$ is $\Theta _{1}\left( \allowbreak \overline{y}\right) +O\left(
\left\| \allowbreak \overline{y}\right\| ^{2}\right) $. Note that we may
extend $f$ to be a basic function on all of $\left( M,\mathcal{F}\right) $
by multiplying by a radial cutoff function and extending by zero. Combining
the two deformations described above, we see that $V_{L}$ may be smoothly
deformed to a $\Gamma $-equivariant transformation of the form $\Theta _{1}$
in such a way that the index is unchanged throughout the deformation.

The argument that follows is somewhat similar to that found in \cite{ESW}.
Let \linebreak $\left\{ Y_{t}:N_{x}L\rightarrow N_{x}L\right\} $ be a smooth,
one-parameter family of $\Gamma $-equivariant transformations constructed
above such that $Y_{0}=V_{L}$ and $Y_{1}\ $is the linear part of $
\sum_{i=1}^{\bar{q}}\pm y^{i}\frac{\partial }{\partial y_{i}}$. Thus each $
Y_{t}\left( \allowbreak \overline{y}\right) $ is a vector field that is
well-defined in a sufficiently small tubular neighborhood of the leaf
closure $L$. Next, let $Z_{t}\left( \allowbreak \overline{y}\right) =h\left(
t\right) \left( \nabla f\left( \allowbreak \overline{y}\right) -Y_{1}\left(
\allowbreak \overline{y}\right) \right) $,with $Y_{1}$ and $\nabla f$ as
above and $h:\Bbb{R\rightarrow }\Bbb{R}$ a smooth positive function such
that $h\left( t\right) =0$ for $t\leq 0$ and $h\left( t\right) =1$ for $
t\geq 1$. Because of the remarks above, $Z_{t}\left( \allowbreak \overline{y}
\right) =O\left( \left\| \allowbreak \overline{y}\right\| ^{2}\right) $ for
all $t$. Let $\eta ^{1}$ be a radial (and therefore basic) cutoff function
that is $1$ in an $r_{1}$-neighborhood of $L$ and $0$ outside an $2r_{1}$
-neighborhood of $L$ ($r_{1}$ will be chosen shortly). Then for $r_{1}$
sufficiently small, $\eta ^{1}\left( \allowbreak \overline{y}\right) \left(
Y_{t}\left( \allowbreak \overline{y}\right) +Z_{t}\left( \allowbreak 
\overline{y}\right) \right) $ is a well-defined basic vector field on $
\left( M,\mathcal{F}\right) $. Let the basic vector field $X_{t}^{1}$ be
defined by 
\begin{equation*}
X_{t}^{1}=\eta ^{1}\left( Y_{t}+Z_{t}\right) +\left( 1-\eta ^{1}\right) V,
\end{equation*}
so that 
\begin{equation*}
X_{t}^{1}\left( \allowbreak \overline{y}\right) =\eta ^{1}\left( \allowbreak 
\overline{y}\right) \left( Y_{t}-Y_{0}\right) \left( \allowbreak \overline{y}%
\right) +Y_{0}\left( y\right) +O\left( \left\| \allowbreak \overline{y}
\right\| ^{2}\right) .
\end{equation*}
Observe that $X_{t}^{1}$ is a smooth, one-parameter family of basic vector
fields that agree with $V$ outside a $2r_{1}$-neighborhood of $L$, and $
X_{1}^{1}=\nabla f$ inside a $r_{1}$-neighborhood of $L$. Let $
m_{1}=\inf_{\left\| \allowbreak \overline{y}\right\| =1,0\leq t\leq
1}\left\| Y_{t}\left( \allowbreak \overline{y}\right) \right\| ,$ and let $
m_{2}=\sup_{\left\| \allowbreak \overline{y}\right\| =1,0\leq t\leq
1}\left\| \frac{d}{dt}Y_{t}\left( \allowbreak \overline{y}\right) \right\| $
. Observe that $0<m_{1}<\infty $ and $0\leq m_{2}<\infty $. Therefore, if $
t_{1}<\frac{m_{1}}{m_{2}}\leq \infty $, then 
\begin{eqnarray*}
\left\| \eta ^{1}\left( \allowbreak \overline{y}\right) \left(
Y_{t}-Y_{0}\right) \left( \allowbreak \overline{y}\right) +Y_{0}\left(
\allowbreak \overline{y}\right) \right\| &\geq &\left\| Y_{0}\left(
\allowbreak \overline{y}\right) \right\| -\left\| \left( Y_{t}-Y_{0}\right)
\left( \allowbreak \overline{y}\right) \right\| \\
&\geq &\left\| Y_{0}\left( \allowbreak \overline{y}\right) \right\| -t\,m_{2}
\\
&>&0
\end{eqnarray*}
for $t\leq t_{1},$ and hence $\eta ^{1}\left( \allowbreak \overline{y}
\right) \left( Y_{t}-Y_{0}\right) \left( \allowbreak \overline{y}\right)
+Y_{0}\left( \allowbreak \overline{y}\right) $ is invertible. Choose $
r_{1}>0 $ so small that $O\left( \left\| \allowbreak \overline{y}\right\|
^{2}\right) <\eta ^{1}\left( \allowbreak \overline{y}\right) \left(
Y_{t}-Y_{0}\right) \left( \allowbreak \overline{y}\right) +Y_{0}\left(
\allowbreak \overline{y}\right) $ for $\left\| \allowbreak \overline{y}
\right\| <2r_{1}$; then the leaf closure $L$ is the only critical leaf
closure of the basic vector field $X_{t}^{1}\left( \allowbreak \overline{y}
\right) $ in the $2r_{1}$-neighborhood of $L$. We continue by defining for $
t_{k}\leq t\leq t_{k+1}:=t_{k}+t_{1}$ 
\begin{eqnarray*}
X_{t}^{k} &=&\eta ^{k}\left( Y_{t}+Z_{t}\right) +\left( 1-\eta ^{k}\right)
X_{t}^{k-1} \\
&=&\eta ^{k}\left( \allowbreak \overline{y}\right) \left(
Y_{t}-Y_{t_{k}}\right) \left( \allowbreak \overline{y}\right)
+Y_{t_{k}}\left( \allowbreak \overline{y}\right) +O\left( \left\|
\allowbreak \overline{y}\right\| ^{2}\right) ,
\end{eqnarray*}
where $\eta ^{k}$ is $1$ in an $r_{k}$-neighborhood of $L$ and $0$ outside
an $2r_{k}$-neighborhood of $L$ (where $0<r_{k}\leq r_{k-1}$ is chosen as
above). Eventually, we will find a $k$ such that $X_{t}^{k}=V$ outside a $
2r_{1}$ -neighborhood of $L$ and is $X_{t}^{k}=Y_{t}+Z_{t}$ inside a $r_{k}$
-neighborhood of $L$. Furthermore, $L$ is the only critical leaf closure of $
X_{t}^{k}$ in a $2r_{1}$ -neighborhood of $L$, and the index of $X_{t}^{k}$
at $L$ is the same as the index of $V$ at $L$ for $0\leq t\leq 1$. Finally, $
X_{1}^{k}=\nabla f$ inside a $r_{k}$-neighborhood of $L$.

Observe that if we let $2r_{1}\leq \delta $, $\widetilde{\delta }=r_{k}$,
and $\widetilde{V}=X_{1}^{k}$, then the first four properties are satisfied.
Moreover, since the orientation line bundles of $Y_{t}$ and $X_{t}^{k}$ at $
L $ are constant in $t$, the last property is satisfied as well.
\end{proof}

We now proceed with a modified version of Witten's deformation of the de
Rham complex (see \cite{W} and \cite{Ro}). Let $V$ be a basic normal vector
field, and let $i\left( V\right) :\Omega _{B}^{*}\left( M\right) \rightarrow
\Omega _{B}^{*}\left( M\right) $ denote interior multiplication with $V$.
For a given $s>0$, let $d_{B,s}=d_{B}+s\,i\left( V\right) :\Omega
_{B}^{*}\left( M\right) \rightarrow \Omega _{B}^{*}\left( M\right) $. Note
that $d_{B,s}$ is the restriction of the differential operator $
d_{s}=d+s\,i\left( V\right) :\Omega ^{*}\left( M\right) \rightarrow \Omega
^{*}\left( M\right) $. The formal adjoint of $d_{B,s}$ is $\delta
_{B,s}=\delta _{B}+sV^{\flat }\wedge :\Omega _{B}^{*}\left( M\right)
\rightarrow \Omega _{B}^{*}\left( M\right) $, where $\delta _{B}$ is the
basic adjoint of $d_{B}$ and $V^{\flat }$ is the basic one-form $
\left\langle V,\cdot \right\rangle $. From \cite{PR}, we know that $\delta
_{B}$ is the restriction of the differential operator $\delta +\varepsilon
:\Omega ^{*}\left( M\right) \rightarrow \Omega ^{*}\left( M\right) $ to
basic forms, where $\varepsilon: \Omega^{i}\left( M\right) \rightarrow
\Omega ^{i-1}\left( M\right)$ is a zeroth-order operator that involves mean
curvature and Rummler's formula. The operator $\varepsilon $ has the
additional property that $P\varepsilon P=0$, where $P:L^{2}\left( \Omega
^{*}\left( M\right) \right) \rightarrow L^{2}\left( \Omega _{B}^{*}\left(
M\right) \right) $ is the orthogonal projection. We define 
\begin{equation*}
D_{B,s}=d_{B,s}+\delta _{B,s}:\Omega _{B}^{*}\left( M\right) \rightarrow
\Omega _{B}^{*}\left( M\right) .
\end{equation*}

Letting $D_{B}=d_{B}+\delta _{B}$ and $H=i\left( V\right) +V^{\flat }\wedge $
, observe that 
\begin{eqnarray*}
D_{B,s}^{2} &=&\left( D_{B}+sH\right) ^{2} \\
&=&D_{B}^{2}+s\left( HD_{B}+D_{B}H\right) +s^{2}H^{2}.
\end{eqnarray*}
The operator $H^{2}=\left( i\left( V\right) +V^{\flat }\wedge \right) ^{2}$
acts by multiplication by the basic function $\left\| V\right\| ^{2}$. A
simple calculation shows $Z^{\prime }:=H(d+\delta +\varepsilon )+(d+\delta
+\varepsilon )H$ satisfies 
\begin{eqnarray}
Z^{\prime } &=&\mathcal{L}_{V}+\left( \mathcal{L}_{V}\right) ^{*,B}+
\left(d\left(V^\flat\right)\right)\wedge +
\left(\left(d\left(V^\flat\right)\right)\wedge \right) ^{*,B}  \notag \\
&=&\mathcal{L}_{V}+\left( \mathcal{L}_{V}\right) ^{*}+Z+
\left(d\left(V^\flat\right)\right)\wedge +
\left(\left((d\left(V^\flat\right)\right)\wedge \right) ^{*}  \label{Zprime}
\end{eqnarray}
where $\mathcal{L}_{V}=i\left( V\right) \,d+d\,i\left( V\right) $ denotes
the Lie derivative in the $V$ direction, the superscript $*,B$ denotes the
adjoint restricted to basic forms, and $Z:=\varepsilon \circ V^{\flat
}\wedge +V^{\flat }\wedge \circ \varepsilon $ is a zeroth order operator.
One can show using the Leibniz rule that $\mathcal{L}_{V}+\left( \mathcal{L}
_{V}\right) ^{*}$ commutes with multiplication by a function, so this
operator is also an operator of order zero. Thus $Z^{\prime }$ is a
differential operator of order zero, and it agrees with $HD_{B}+D_{B}H$ on
basic forms. Also note that $Z^{\prime }$ maps odd forms to odd forms and
even forms to even forms.

Observe that $D_{B}^{2}=\Delta _{B}^{j}=d_{B}\delta _{B}+\delta _{B}d_{B}$
on basic $j$-forms, the basic Laplacian. By the results of \cite{PR}, this
operator is essentially self-adjoint and has eigenvalues $0\leq \lambda
_{1}^{B,j}\leq \lambda _{2}^{B,j}\leq \lambda _{3}^{B,j}\leq ...$ with the
property that $\lambda _{k}^{B,j}\geq Ck^{2/n}$ for some positive constant $
C $ and sufficiently large $k$ (see \cite{R} for more precise asymptotics).
Furthermore, the basic Hodge decomposition theorem (see \cite{PR}, \cite{KT}
) states that $\ker \Delta _{B}$ is finite-dimensional and that the space of
basic $j$-forms decomposes orthogonally as $\mathrm{im}\,d_{B}\oplus \mathrm{%
\ \ im}\,\delta _{B}\oplus \ker \Delta _{B}^{j}$. Letting $\Delta _{B}^{j}$
denote the Laplacian on $j$-forms, we have $\ker \Delta _{B}^{j}\cong
H_{B}^{j}\left( M,\mathcal{F}\right) $. This implies 
\begin{eqnarray*}
\allowbreak \chi _{B}\left( M,\mathcal{F}\right) &=&\sum_{j=0}^{q}\left(
-1\right) ^{j}\dim H_{B}^{j}\left( M,\mathcal{F}\right) \\
&=&\sum_{j=0}^{q}\left( -1\right) ^{j}\dim \ker \Delta _{B}^{j} \\
&=&\mathrm{index}\left( D_{B}:\Omega _{B}^{\mathrm{even}}\left( M\right)
\rightarrow \Omega _{B}^{\mathrm{odd}}\left( M\right) \right) .
\end{eqnarray*}

Next, let $K_{B}^{j}\left( t,x,y\right) $ denote the basic heat kernel on $j$
-forms (see \cite{KT},\cite{PR},\cite{R}), which is the fundamental solution
of the basic heat equation. More specifically, consider the bundle $
E\rightarrow \left[ 0,\infty \right) \times M\times M$ where $E_{\left(
t,x,y\right) }=\mathrm{Hom}\left( \bigwedge^{j}\left( N_{y}\mathcal{F}
\right) ^{*},\bigwedge^{j}\left( N_{x}\mathcal{F}\right) ^{*}\right) $. The
basic heat kernel is a section of $E$ that is basic in $x$ and $y$ and
satisfies 
\begin{eqnarray*}
\left( \frac{\partial }{\partial t}+\Delta _{B,x}^{j}\right) K_{B}^{j}\left(
t,x,y\right) &=&0\text{ for }t>0; \\
\lim_{t\rightarrow 0^{+}}\int_{M_{y}}K_{B}^{j}\left( t,x,y\right) \beta
\left( y\right) \,\omega _{\mathrm{vol}}\left( y\right) &=&\beta \left(
x\right) \text{ for all basic }j\text{-forms }\beta .
\end{eqnarray*}
In \cite[Theorem 3.5]{PR}, the authors show that $K_{B}^{j}\left(
t,x,y\right) $ exists, is smooth in $x,y,$ and $t$, is unique, and satisfies 
\begin{equation*}
K_{B}^{j}\left( t,x,y\right) =P_{x}P_{y}\tilde{K}^{j}\left( t,x,y\right)
=\sum_{k=1}^{\infty }e^{-\lambda _{k}^{B,j}t}\alpha _{k}\left( x\right)
\otimes \alpha _{k}^{*}\left( y\right) ,
\end{equation*}
where the set of $j$-forms $\left\{ \alpha _{1},\alpha _{2},...\right\} $ is
an orthonormal basis of basic eigenforms corresponding to the eigenvalues $
\left\{ \lambda _{1}^{B,j},\lambda _{2}^{B,j},...\right\} $, and $\tilde{K}
^{j}\left( t,x,y\right) $ is the heat kernel corresponding to the strongly
elliptic operator $\Delta ^{j}+\delta \varepsilon ^{*}+\varepsilon
^{*}\delta $. The operator $\Delta ^{j}$ is the ordinary Laplacian on 
$j$-forms. Note that this sum is finite if the leaves are dense. The map $P_{x}$
is the orthogonal projection from $L^{2}\left( \Omega ^{*}\left( M\right)
\right) $ to $L^{2}\left( \Omega _{B}^{*}\left( M\right) \right) $ in the $x$
-variable. The map $P_{y}$ is the induced basic projection on duals of forms
in the $y$-variable. In \cite{PR}, the authors show that this basic
projection $P$ (or $P_{x}$ or $P_{y}$) maps smooth forms to smooth basic
forms, and the results in \cite{PR} also imply that the map $\beta
\longmapsto P\beta $ is smooth. The main other fact used in the proofs of
the results concerning the basic Laplacian and the basic heat kernel is that 
$\Delta _{B}^{j}$ is the restriction of $\Delta ^{j}+\varepsilon
d+d\varepsilon $ , which is a strongly elliptic operator defined on all
forms.

The standard heat kernel approach to index calculations carries over to the
basic case. Since $D_{B}$ maps an eigenspace of $\Delta _{B}$ in $\Omega
_{B}^{\mathrm{even}}\left( M\right) $ isomorphically onto an eigenspace in $
\Omega _{B}^{\mathrm{odd}}\left( M\right) $, we have 
\begin{eqnarray*}
\allowbreak \chi _{B}\left( M,\mathcal{F}\right) &=&\mathrm{index}\left(
D_{B}:\Omega _{B}^{\mathrm{even}}\left( M\right) \rightarrow \Omega _{B}^{%
\mathrm{odd}}\left( M\right) \right) \\
&=&\sum_{j=0}^{q}\left( -1\right) ^{j}\mathrm{tr}\left( e^{-t\Delta
_{B}^{j}}\right) \\
&=&\sum_{j=0}^{q}\left( -1\right) ^{j}\int_{M}\mathrm{tr}K_{B}^{j}\left(
t,x,x\right) \,\,\omega _{\mathrm{vol}}\left( x\right) \\
&=&\int_{M}\mathrm{tr}K_{B}^{\mathrm{even}}\left( t,x,x\right) \,\,\omega _{%
\mathrm{vol}}\left( x\right) -\int_{M}\mathrm{tr}K_{B}^{\mathrm{odd}}\left(
t,x,x\right) \,\,\omega _{\mathrm{vol}}\left( x\right) .
\end{eqnarray*}

A similar analysis may be applied to $D_{B,s}=D_{B}+sH$. The operator $
D_{B,s}^{2}=D_{B}^{2}+s\left( HD_{B}+D_{B}H\right) +s^{2}\left\| V\right\|
^{2}$ is the restriction of the strongly elliptic operator 
\begin{equation}
\Delta _{s}^{\prime }=\Delta +\varepsilon d+d\varepsilon +sZ^{\prime
}+s^{2}\left\| V\right\| ^{2}  \label{restrictD^2}
\end{equation}
to basic forms, and the terms involving $s$ are zeroth order. By a proof
similar to that in \cite{PR}, this operator is essentially self-adjoint and
has eigenvalues that grow at a rate similar to those of $\Delta _{B}$. The
basic heat kernel $K_{B,s}^{\mathrm{even}}\left( t,x,y\right) $
corresponding to $D_{B,s}^{2}$ exists, is smooth in $x,y,$ and $t$, is
unique, and satisfies 
\begin{equation*}
K_{B}^{\mathrm{even}}\left( t,x,y\right) =P_{x}P_{y}\tilde{K}_{s}^{\mathrm{\
even}}\left( t,x,y\right) =\sum_{k=1}^{\infty }e^{-t\,\lambda _{k}^{B,s,%
\mathrm{even}}}\alpha _{k}^{s}\left( x\right) \otimes \left( \alpha
_{k}^{s}\left( y\right) \right) ^{*},
\end{equation*}
where the set $\left\{ \alpha _{1}^{s},\alpha _{2}^{s},...\right\} $ is an
orthonormal basis of even basic eigenforms corresponding to the eigenvalues $
\left\{ \lambda _{1}^{B,s,\mathrm{even}},\lambda _{2}^{B,s,\mathrm{even}
},...\right\} $ of $D_{B,s}^{2}$, and $\tilde{K}_{s}^{\mathrm{even}}\left(
t,x,y\right) $ is the heat kernel corresponding to the strongly elliptic
operator $\Delta ^{\mathrm{even}}+\delta \varepsilon ^{*}+\varepsilon
^{*}\delta +s\left( Z^{\prime }\right) ^{*}+s^{2}\left\| V\right\| ^{2}$.
Similar results are true for $D_{B,s}^{2}$ on odd forms. Again, we have that 
$D_{B,s}$ maps an eigenspace of $D_{B,s}^{2}$ in $\Omega _{B}^{\mathrm{even}
}\left( M\right) $ isomorphically onto an eigenspace in $\Omega _{B}^{%
\mathrm{odd}}\left( M\right) $, so that 
\begin{eqnarray}
&&\mathrm{index}\left( D_{B,s}:\Omega _{B}^{\mathrm{even}}\left( M\right)
\rightarrow \Omega _{B}^{\mathrm{odd}}\left( M\right) \right)  \notag \\
&=&\int_{M}\mathrm{tr}K_{B,s}^{\mathrm{even}}\left( t,x,x\right) \,\,\omega
_{\mathrm{vol}}\left( x\right) -\int_{M}\mathrm{tr}K_{B,s}^{\mathrm{odd}
}\left( t,x,x\right) \,\,\omega _{\mathrm{vol}}\left( x\right) .
\label{indexkernel}
\end{eqnarray}

The above results about $D_{B,s}^{2}$ imply results about the first order
operator $D_{B,s}$. The spectrum of $D_{B,s}^{2}$ is nonnegative; let $
E_{\lambda ^{2}}$ be an eigenspace of $D_{B,s}^{2}$ with eigenvalue $\lambda
^{2}$. Then $\left( D_{B,s}+\lambda \right) E_{\lambda ^{2}}$ and $\left(
D_{B,s}-\lambda \right) E_{\lambda ^{2}}$ are finite-dimensional orthogonal
eigenspaces of $D_{B,s}$ with eigenvalues $\lambda $ and $-\lambda $,
respectively. Since $D_{B,s}$ is the restriction of the first order elliptic
operator 
\begin{equation*}
D_{s}=d+\delta +\varepsilon +s\left( i\left( V\right) +V^{\flat }\wedge
\right) ,
\end{equation*}
the subspaces $\left( D_{B,s}\pm \lambda \right) E_{\lambda ^{2}}$ are
spanned by smooth eigenforms. Moreover, there is an orthonormal basis of $
L^{2}\left( \Omega _{B}^{*}\left( M\right) \right) $ consisting of such
smooth eigenforms.

The operators $\frac{\partial }{\partial t}-iD_{s}^{*}$ and $\frac{\partial 
}{\partial t}+iD_{s}^{*}$ are strongly hyperbolic, so our initial value
problem has a unique solution (see \cite[sections 69-74]{D}, 
\cite[chapters IV-V]{La}, \cite[chapter 6]{Mi}, \cite[section 6.5]{Ta}).
Thus, we define the function $e^{itD_{s}^{*}}u$ to be the unique solution to
the initial value problem (the \emph{generalized traveling wave equation}) 
\begin{eqnarray*}
\left( \frac{\partial }{\partial t}-iD_{s}^{*}\right) f\left( t,x\right)
&=&0\,\text{\thinspace \thinspace for all }t>0,x\in M \\
f\left( 0,x\right) &=&u\left( x\right) \,\text{\thinspace for all }x\in M.
\end{eqnarray*}
Since $D_{s}$ maps basic forms to basic forms, $PD_{s}P=D_{s}P$, and $
PD_{s}^{*}P=PD_{s}^{*}$; the form $e^{itD_{s}}u$ is basic for a given basic
form $u$. Moreover, 
\begin{equation*}
\frac{\partial }{\partial t}\left( e^{itD_{s}}u\right) =iD_{s}\left(
e^{itD_{s}}u\right) =iD_{B,s}\left( e^{itD_{s}}u\right) ,
\end{equation*}
so that $e^{itD_{B,s}}u:=e^{itD_{s}}u$ is a solution to the \emph{basic
traveling wave equation} with the appropriate initial condition: 
\begin{eqnarray}
\left( \frac{\partial }{\partial t}-iD_{B,s}\right) \beta \left( t,x\right)
&=&0\,\text{\thinspace \thinspace for all }t>0,x\in M  \notag \\
\beta \left( 0,x\right) &=&u\left( x\right) \in \Omega _{B}^{*}\left(
M\right) .  \label{basictravwave}
\end{eqnarray}
Then $e^{itD_{B,s}}u$ and $e^{-itD_{B,s}}u$ are both solutions of the basic
wave equation 
\begin{equation}
\left( \frac{\partial ^{2}}{\partial t^{2}}+D_{B,s}^{2}\right) \beta \left(
t,x\right) =0\,.  \label{basicwave}
\end{equation}
Because of the results concerning $\frac{\partial }{\partial t}-iD_{s}$, the
solutions to these wave equations exist and are unique with respect to the
appropriate initial conditions (both $\beta \left( 0,x\right) $ and $\frac{%
\partial }{\partial t}\beta \left( 0,x\right) $ must be specified for the
basic wave equation).

We now prove some analytic results about the operators $D_{B,s}$ and similar
operators; we include proofs only where the standard proofs do not translate
directly to the basic case. Let $L^{p}=L^{p}\left( \Omega _{B}^{*}\left(
M\right) \right) $ denote the $L^{p}-$norm closure of the space of smooth
basic forms, and let $W^{k}=W^{k}\left( \Omega _{B}^{*}\left( M\right)
\right) $ denote the closure of the space of these basic forms under the
Sobolev $\left( k,2\right) $-norm. Let $\left\| \cdot \right\| _{p}$ and $
\left\| \cdot \right\| _{k,2}$ denote the norms on these spaces.

\begin{lemma}
\label{elliptic1}Let $D$ be a strongly elliptic, first order operator on $
\Omega ^{\ast }\left( M\right) $ that restricts to an operator on $\Omega
_{B}^{\ast }\left( M\right) $. Suppose that the restriction of $D$ is
formally self-adjoint on $\Omega _{B}^{\ast }\left( M\right) $. Then, for
some $c>0,$ $\left\| \left( D^{2}+1\right) \alpha \right\| _{k,2}\geq
c\left\| \alpha \right\| _{k+1,2}$ for every basic form $\alpha $.
\end{lemma}

\begin{proof}
We use induction on $k$. For $k=0$, we clearly have that 
\begin{eqnarray*}
\left\| \left( D^{2}+1\right) \alpha \right\| _{0,2} &=&\left( \left\langle
D^{2}\alpha ,D^{2}\alpha \right\rangle +2\left\langle D\alpha ,D\alpha
\right\rangle +\left\langle \alpha ,\alpha \right\rangle \right) ^{1/2} \\
&\geq &\left\| \alpha \right\| _{0,2}.
\end{eqnarray*}
Similarly, $\left\| \left( D^{2}+1\right) \alpha \right\| _{0,2}\geq \sqrt{2}
\left\| D\alpha \right\| _{0,2}$. By ellipticity, there exist $c_{1}>0$ and $
c_{2}\geq 0$ (independent of $\alpha $) such that $\sqrt{2}\left\| D\alpha
\right\| _{0,2}\geq c_{1}\left\| \alpha \right\| _{1,2}-c_{2}\left\| \alpha
\right\| _{0,2}$. Thus, 
\begin{equation*}
\left\| \left( D^{2}+1\right) \alpha \right\| _{0,2}\geq \left( \frac{c_{1}}{
1+c_{2}}\right) \left\| \alpha \right\| _{1,2}.
\end{equation*}

Assume that the conclusion is true for $0\leq k\leq m$. Then 
\begin{eqnarray*}
\left\| \left( D^{2}+1\right) \alpha \right\| _{m+1,2} &\geq &c_{1}\left\|
\left( D^{2}+1\right) ^{2}\alpha \right\| _{m-1,2}\text{ by definition of
Sobolev norm} \\
&\geq &c_{2}\left\| \left( D^{2}+1\right) \alpha \right\| _{m,2}\text{ by
the induction hypothesis} \\
&\geq &c_{3}\left\| \alpha \right\| _{m+2,2}-c_{4}\left\| \alpha \right\|
_{m+1,2}\text{ by ellipticity,}
\end{eqnarray*}
for positive constants $c_{1}$, $c_{2}$, and $c_{3}$ and for $c_{4}\geq 0$,
all independent of $\alpha $. Also, it follows from our induction hypothesis
that $\left\| \left( D^{2}+1\right) \alpha \right\| _{m,2}\geq c_{5}\left\|
\alpha \right\| _{m+1,2}$ for some $c_{5}>0$. Finally, by the definition of
the Sobolev norm, there exists a positive constant $c_{6}$ such that 
\begin{eqnarray*}
\left\| \left( D^{2}+1\right) \alpha \right\| _{m+1,2} &\geq &c_{6}\left\|
\left( D^{2}+1\right) \alpha \right\| _{m+1,2}+\frac{c_{6}c_{4}}{c_{5}}%
\left\| \left( D^{2}+1\right) \alpha \right\| _{m,2} \\
&\geq &c_{6}c_{3}\left\| \alpha \right\| _{m+2,2}\text{,}
\end{eqnarray*}
whence the result follows.
\end{proof}

\begin{lemma}
\label{elliptic2}Let $D$ be as in Lemma \ref{elliptic1}. Then there exists a
constant $c>0$ such that $\left\| \left( D^{2}+1\right) \alpha \right\|
_{k,2}\geq c\left\| \alpha \right\| _{k+2,2}$ for every basic form $\alpha $.
\end{lemma}

\begin{proof}
By ellipticity, there exist $c_{1}>0$ and $c_{2}\geq 0$ such that $\left\|
\left( D^{2}+1\right) \alpha \right\| _{k,2}\geq c_{1}\left\| \alpha
\right\| _{k+2,2}-c_{2}\left\| \alpha \right\| _{k+1,2}$. By Lemma \ref
{elliptic1}, there exists $c_{3}>0$ such that $\left\| \left( D^{2}+1\right)
\alpha \right\| _{k,2}\geq c_{3}\left\| \alpha \right\| _{k+1,2}$. Thus 
\begin{equation*}
\left\| \left( D^{2}+1\right) \alpha \right\| _{k,2}\geq \left( \frac{c_{1}}{%
1+\frac{c_{2}}{c_{3}}}\right) \left\| \alpha \right\| _{k+2,2}.
\end{equation*}
\end{proof}

\begin{remark}
\label{ellipticremark}Observe that if the coefficients of the operator $D$
depend polynomially on a parameter $s$, then the constants in the
inequalities of the lemmas above can be chosen to depend polynomially in $s$
. Also, the results of Lemma \ref{elliptic1} and Lemma \ref{elliptic2}
extend to much more general situations; the operator $D$ may be any strongly
elliptic, first-order operator acting on smooth sections of a vector bundle
over a compact manifold such that the restriction of $D$ to a subspace
consisting of smooth sections is formally self-adjoint.
\end{remark}

\begin{lemma}
(Basic Sobolev Embedding Theorem) \label{sobolev}If $a>\frac{n}{2}$, then $
W^{a}\subset L^{\infty }$. Also, $W^{a+k}\subset C^{k}$ for $k>0$.
\end{lemma}

\begin{proof}
The standard proof (see \cite{Ro}) works for the basic case.
\end{proof}

\begin{lemma}
\label{dualL^1}Define $L_{\beta }\left( \alpha \right) =\int_{M}\left(
\alpha ,\beta \right) \,\omega _{\mathrm{vol}}.$ Then the map $\beta \mapsto
L_{\beta }$ defines a norm-preserving injection of $L^{1}\left( \Omega
_{B}^{\ast }\left( M\right) \right) $ into $\left( L^{\infty }\left( \Omega
_{B}^{\ast }\left( M\right) \right) \right) ^{\ast }$.
\end{lemma}

\begin{proof}
The standard proof (see \cite{Ro}) works for the basic case.
\end{proof}

Recall that a subset of a foliation is called \emph{saturated} if it is a
union of leaves.

\begin{lemma}
\label{saturatedFubini}Let $U$ be a saturated open set in a Riemannian
foliation $\left( M,\mathcal{F}\right) $. Let $P:L^{2}\left( \Omega ^{\ast
}\left( M\right) \right) \rightarrow L^{2}\left( \Omega _{B}^{\ast }\left(
M\right) \right) $ be the orthogonal projection. Then for every $f\in
C^{\infty }\left( M\right) ,$ 
\begin{equation*}
\int_{U}f\,\omega _{\mathrm{vol}}=\int_{U}Pf\,\omega _{\mathrm{vol}}.
\end{equation*}
\end{lemma}

\begin{proof}
Let $\left\{ \eta _{\varepsilon }\right\} $ be a family of smooth basic
functions approximating (in $L^{2}$) the characteristic function of $U$ as $
\varepsilon \rightarrow 0$ (see the results of \cite{PR} and \cite{LeR}).
Then $\int_{M}f\eta _{\varepsilon }\,\omega _{\mathrm{vol}}=\int_{M}P\left(
f\eta _{\varepsilon }\right) \,\omega _{\mathrm{vol}}=\int_{M}\left(
Pf\right) \eta _{\varepsilon }\,\omega _{\mathrm{vol}}$ since $\eta
_{\varepsilon }$ is basic (\cite{PR}). We then apply the dominated
convergence theorem.
\end{proof}

Let $\left( \cdot ,\cdot \right) $ denote the pointwise inner product of
forms. The following two results apply to $D_{B,s}$, since $D_{B,s}$ is the
restriction of $D_{s}=d+\delta +\varepsilon +s\left( i\left( V\right)
+V^{\flat }\wedge \right) $ to basic forms and since $P\varepsilon P=0$ (see 
\cite{PR}).

\begin{lemma}
\label{EnergyEst}Let $D$ be a first order operator on $\Omega ^{\ast }\left(
M\right) $ of the form $D=d+\delta +Z_{1}+Z_{2}$, where

\begin{enumerate}
\item  $D$ restricts to an operator on $\Omega _{B}^{\ast }\left( M\right) $;

\item  The operator $Z_{1}$ is zeroth order and satisfies $PZ_{1}P=0$, where 
$P$ is the orthogonal projection from $L^{2}\left( \Omega ^{\ast }\left(
M\right) \right) $ to $L^{2}\left( \Omega _{B}^{\ast }\left( M\right)
\right) $;

\item  The operator $Z_{2}$ is zeroth order and is formally self-adjoint
with respect to the pointwise inner product of forms.
\end{enumerate}

Let $B\left( L,R\right) $ denote the set of points of distance less than $R$
from a fixed leaf closure $L$; assume that $R$ is sufficiently small to
avoid the cut locus of $L$. Let $\beta $ be a basic form,
and let $\beta _{t}=e^{itD}\beta $ be the corresponding solution to the
basic traveling wave equation as in (\ref{basictravwave}). Then, for $0\leq
t<R$, the function 
\begin{equation*}
f\left( t\right) =\int_{B\left( L,R-t\right) }\left( \beta _{t},\beta
_{t}\right) \,\omega _{\mathrm{vol}}
\end{equation*}
is a decreasing function of $t$.
\end{lemma}

\begin{proof}
We have 
\begin{equation*}
\frac{df}{dt}=\int_{B\left( L,R-t\right) }\left( \left( iD\beta _{t},\beta
_{t}\right) +\left( \beta _{t},iD\beta _{t}\right) \right) \,-\int_{S\left(
L,R-t\right) }\left( \beta _{t},\beta _{t}\right) \,,
\end{equation*}
where $S\left( L,r\right) $ is the set of points of distance $R$ from the
leaf closure $L$. Observe that 
\begin{eqnarray*}
\left( iD\beta _{t},\beta _{t}\right) +\left( \beta _{t},iD\beta _{t}\right)
\, &=&\left( \left( i\left( d+\delta +Z_{1}+Z_{2}\right) \beta _{t},\beta
_{t}\right) +\left( \beta _{t},i\left( d+\delta +Z_{1}+Z_{2}\right) \beta
_{t}\right) \right) \\
&=&i\left[ \left( \left( d+\delta \right) \beta _{t},\beta _{t}\right)
-\left( \beta _{t},\left( d+\delta \right) \beta _{t}\right) \right]
\\
&\ &+i\left[ \left( Z_{1}\beta _{t},\beta _{t}\right) -\left( \beta
_{t},Z_{1}\beta _{t}\right) \right] ,
\end{eqnarray*}
since $Z_{2}$ is formally self-adjoint with respect to $\left( \cdot ,\cdot
\right) $. By \cite[proof of Proposition 2.9]{Ro}, 
\begin{equation*}
i\left[ \left( \left( d+\delta \right) \beta _{t},\beta _{t}\right) -\left(
\beta _{t},\left( d+\delta \right) \beta _{t}\right) \right] =i\delta \omega
_{t},
\end{equation*}
where $\omega $ is the one-form defined by $\omega _{t}\left( X\right)
=-\left( X\mathbf{.}\beta _{t},\beta _{t}\right) :=-\left( \left( X^{\flat
}\wedge -i\left( X\right) \right) \beta _{t},\beta _{t}\right) $. Then 
\begin{eqnarray*}
\int_{B\left( L,R-t\right) }\left( \left( iD\beta _{t},\beta _{t}\right)
+\left( \beta _{t},iD\beta _{t}\right) \right) &=&\int_{B\left( L,R-t\right)
}i\delta \omega _{t}+i\left[ \left( Z_{1}\beta _{t},\beta _{t}\right)
-\left( \beta _{t},Z_{1}\beta _{t}\right) \right] \\
&=&\int_{B\left( L,R-t\right) }i\delta \omega _{t}+iP\left[ \left(
Z_{1}P\beta _{t},\beta _{t}\right) -\left( \beta _{t},Z_{1}P\beta
_{t}\right) \right] ,
\end{eqnarray*}
which follows from Lemma \ref{saturatedFubini} and the fact that $\beta _{t}$
is basic. Furthermore, by the results of \cite{PR}, $P\left[ \left(
Z_{1}P\beta _{t},\beta _{t}\right) -\left( \beta _{t},Z_{1}P\beta
_{t}\right) \right] =\left( PZ_{1}P\beta _{t},\beta _{t}\right) -\left(
\beta _{t},PZ_{1}P\beta _{t}\right) $. Since $PZ_{1}P=0$ by hypothesis, the
divergence theorem yields 
\begin{eqnarray*}
\frac{df}{dt} &=&i\int_{B\left( L,R-t\right) }\delta \omega
_{t}\,-\int_{S\left( L,R-t\right) }\left( \beta _{t},\beta _{t}\right) \\
&=&i\int_{S\left( L,R-t\right) }\omega _{t}\left( N\right) \,-\int_{S\left(
L,R-t\right) }\left( \beta _{t},\beta _{t}\right) ,
\end{eqnarray*}
where $N$ is the unit vector field normal to $S\left( L,R-t\right) $ with
orientation chosen compatibly with the choice of orientation of $S\left(
L,R-t\right) $. Locally $\left\| \omega _{t}\left( N\right) \right\|
^{2}=\left\| \left( N\mathbf{.}\beta _{t},\beta _{t}\right) \right\|
^{2}\leq \left( \beta _{t},\beta _{t}\right) \left( N\mathbf{.}\beta _{t},N 
\mathbf{.}\beta _{t}\right) =\left( \beta _{t},\beta _{t}\right) ^{2}$.
Since $\frac{df}{dt}$ is real, we conclude that 
\begin{equation*}
\frac{df}{dt}\leq 0.
\end{equation*}
\end{proof}

\begin{proposition}
\label{UnitPropSpeed}(Unit Propagation Speed) Let $D$ be as in Lemma \ref
{EnergyEst}.Then for any $\beta \in \Omega _{B}^{\ast }\left( M\right) $,
the support of $e^{itD}\beta $ lies within a distance $\left| t\right| $ of
the support of $\beta $.
\end{proposition}

\begin{proof}
One easily checks that such an operator $D$ is formally self-adjoint on the
space of basic forms, and $e^{itD}e^{\pm isD}\beta =e^{i\left( t\pm s\right)
D}\beta $; thus, it is sufficient to prove the result for small positive $t$
. Since $\beta $ is basic, the support of $\beta $ and its complement are
saturated. Since $\left( M,\mathcal{F}\right) $ is a Riemannian foliation
and since $M$ is compact, for every leaf closure $L$, there exists a tubular
neighborhood and $R>0$ such that for every leaf closure $L_{x}$ in that
tubular neighborhood, the set $B\left( L_{x},R\right) $ of points of
distance less than $R$ from $L_{x}$ misses the focal locus and cut locus of $
L_{x}$. Choose any $x\in M$ that is at a distance $R$ or more from the
support of $\beta $; let $L_{x}$ denote the leaf closure containing $x$.
Since $\left( M,\mathcal{F}\right) $ is Riemannian, the set $B\left(
L_{x},R\right) $ is also disjoint from the support of $\beta $. Then 
\begin{equation*}
\int_{B\left( L_{x},R\right) }\left( \beta ,\beta \right) =0\geq
\int_{B\left( L_{x},R-t\right) }\left( e^{itD}\beta ,e^{itD}\beta \right) 
\text{ for }0<t<R,
\end{equation*}
by Lemma \ref{EnergyEst}. Hence $e^{itD}\beta =0$ at $x$ for $0<t<R$.
\end{proof}

\begin{lemma}
\label{smoothkernel}Let $D$ be as in Lemma \ref{EnergyEst}, and let $\phi : 
\Bbb{R}\rightarrow \Bbb{R}$ be a rapidly decreasing even function. Then $
\phi \left( D\right) $ has a smooth kernel $k\left( x,y\right) $ that is
basic in each factor and satisfies $k\left( x,y\right) =\overline{k\left(
y,x\right) }$.
\end{lemma}

\begin{proof}
The operator $\phi \left( D\right) $ is a bounded, self-adjoint operator on $
L^{2}\left( \Omega _{B}^{\ast }\right) $. By Lemma \ref{elliptic1} and Lemma 
\ref{elliptic2}, $\phi \left( D\right) $ maps $L^{2}\left( \Omega _{B}^{\ast
}\right) $ into $W^{k}$ for every $k>0$. By Lemma \ref{sobolev}, this
implies $\phi \left( D\right) $ maps $L^{2}\left( \Omega _{B}^{\ast }\right) 
$ to smooth basic forms. Using the fact that $D$ is formally self-adjoint on 
$\Omega _{B}^{\ast }\left( M\right) $, the standard proofs are easily
adapted to show the existence and properties of $k\left( x,y\right) $ (see 
\cite[Lemma 5.6]{Ro}).
\end{proof}

The index of a Fredholm operator is a homotopy invariant; we prove a similar
result for the operators $D_{B,s}$ on basic forms. The operators $\left(
\Delta _{s}^{\prime }\right) ^{*}=\Delta ^{j}+\delta \varepsilon
^{*}+\varepsilon ^{*}\delta +s\left( Z^{\prime }\right) ^{*}+s^{2}\left\|
V\right\| ^{2}$ depend smoothly on the parameter $s$, the vector field $V$,
and the metric. Therefore the kernels $\tilde{K}_{s}^{\mathrm{even}}\left(
t,x,y\right) $ and $\tilde{K}_{s}^{\mathrm{odd}}\left( t,x,y\right) $ are
smooth in these parameters as well. Since the map $\beta \longmapsto P\beta $
is continuous, $K_{B,s}^{\mathrm{even}}\left( t,x,y\right) $ and $K_{B,s}^{%
\mathrm{odd}}\left( t,x,y\right) $ are continuous in $s$, $V$, and the
metric. By (\ref{indexkernel}) and the fact that the index is an integer, we
have shown the following:

\begin{proposition}
\label{indexeuler}$\allowbreak \chi _{B}\left( M,\mathcal{F}\right) =\mathrm{%
index}\left( D_{B,s}:\Omega _{B}^{\mathrm{even}}\left( M\right) \rightarrow
\Omega _{B}^{\mathrm{odd}}\left( M\right) \right) $ for all $s\in \Bbb{R}$.
\end{proposition}\qed

Let $\phi $ be a smooth, rapidly decreasing function on $[0,\infty )$ with $
\phi \left( 0\right) =1.$ Then $\phi \left( D_{B}^{2}\right) $ is a trace
class operator. Let 
\begin{equation*}
\mu _{j}=\mathrm{tr}\left( \left. \phi \left( D_{B}^{2}\right) \right|
_{L^{2}\left( \Omega _{B}^{j}\left( M\right) \right) }\right).
\end{equation*}
By adapting the argument of \cite[Chapter 12]{Ro} for the standard case, we
obtain

\begin{theorem}
(Basic Morse Inequalities) \label{basicmorseineq}The numbers $\mu _{j}$ and $
\beta _{j}=\dim H_{B}^{j}\left( M,\mathcal{F}\right) $ satisfy the following
system of inequalities: 
\begin{eqnarray*}
\beta _{0} &\leq &\mu _{0}, \\
\beta _{1}-\beta _{0} &\leq &\mu _{1}-\mu _{0}, \\
\beta _{2}-\beta _{1}+\beta _{0} &\leq &\mu _{2}-\mu _{1}+\mu _{0,}
\end{eqnarray*}
etc., and the equality 
\begin{eqnarray*}
\allowbreak \chi _{B}\left( M,\mathcal{F}\right) &=&\mathrm{index}\left(
D_{B}:\Omega _{B}^{\mathrm{even}}\left( M\right) \rightarrow \Omega _{B}^{%
\mathrm{odd}}\left( M\right) \right) \\
&=&\sum_{j=0}^{q}\left( -1\right) ^{j}\beta _{j}=\sum_{j=0}^{q}\left(
-1\right) ^{j}\mu _{j}.
\end{eqnarray*}
\end{theorem}

\begin{proof}
The proof is identical to the proof of the standard case in \cite[Chapter 12]
{Ro}, replacing $d$ with $d_{B}$ and the standard Hodge theorem with the
basic Hodge theorem.
\end{proof}

Let 
\begin{eqnarray*}
\mu _{s}^{\mathrm{even}} &=&\mathrm{tr}\left( \left. \phi \left(
D_{B,s}^{2}\right) \right| _{L^{2}\left( \Omega _{B}^{\mathrm{even}}\left(
M\right) \right) }\right) \\
\mu _{s}^{\mathrm{odd}} &=&\mathrm{tr}\left( \left. \phi \left(
D_{B,s}^{2}\right) \right| _{L^{2}\left( \Omega _{B}^{\mathrm{odd}}\left(
M\right) \right) }\right) .
\end{eqnarray*}
Combining the proof of the standard Morse inequalities with Proposition~\ref
{indexeuler}, we have the following result.

\begin{proposition}
\allowbreak \label{indepofs}\label{eulermu}$\chi _{B}\left( M,\mathcal{F}
\right) =\mu _{s}^{\mathrm{even}}-\mu _{s}^{\mathrm{odd}}$ for all $s\in 
\Bbb{R}$.
\end{proposition}\qed

Fix a number $\rho >0$, and choose a positive, even Schwarz function $\phi $
with $\phi \left( 0\right) =1$ and such that the Fourier transform $\hat{\phi%
}\left( \xi \right) =\int_{\Bbb{R}}e^{-ix\xi }\phi \left( x\right) \,dx$ is
supported in the interval $\left[ -\rho ,\rho \right] $. Since $\phi $ is
even, $\phi \left( D_{B,s}\right) $ makes sense as a smoothing operator, and
the basic Euler characteristic satisfies Proposition \ref{eulermu} with 
\begin{eqnarray*}
\mu _{s}^{\mathrm{even}} &=&\mathrm{tr}\left( \left. \phi \left(
D_{B,s}\right) \right| _{L^{2}\left( \Omega _{B}^{\mathrm{even}}\left(
M\right) \right) }\right) \text{ and} \\
\mu _{s}^{\mathrm{odd}} &=&\mathrm{tr}\left( \left. \phi \left(
D_{B,s}\right) \right| _{L^{2}\left( \Omega _{B}^{\mathrm{odd}}\left(
M\right) \right) }\right) .
\end{eqnarray*}
Let $\mathrm{Crit}\left( V\right) $ be the (finite) union of critical leaf
closures of $V$ in $M$.

\begin{lemma}
\label{complement}On the complement of a $2\rho $-neighborhood of $\mathrm{%
Crit}\left( V\right) $, the basic kernel of $\phi \left( D_{B,s}\right) $
satisfies $k_{B,s}\left( x,y\right) \rightarrow 0$ uniformly as $
s\rightarrow \infty $.
\end{lemma}

\begin{proof}
This proof is the same as \cite[Lemma 12.10]{Ro}, but we include some
details here because of subtleties. Choose a constant $C$ so that $\left\|
V\left( x\right) \right\| \geq C>0$ for all $x$ in the complement of a $\rho 
$-neighborhood of $\mathrm{Crit}\left( V\right) $. By formula (\ref
{restrictD^2}) and the remarks following that equation, $\left\langle
D_{B,s}^{2}\beta ,\beta \right\rangle \geq \frac{C^{2}s^{2}}{2}\left\langle
\beta ,\beta \right\rangle $ for every $\beta \in \Omega _{B}^{*}\left(
M\right) $ that is supported on the complement of such a neighborhood and
for sufficiently large $s$. Let $\mathcal{H}$ denote the Hilbert space of $
L^{2}$ basic forms that vanish on a $\rho $-neighborhood of $\mathrm{Crit}
\left( V\right) $. Then $D_{B,s}^{2}$ is a positive symmetric operator on a
dense subset of $\mathcal{H}$, so it extends to a self-adjoint operator $A$
on $\mathcal{H}$ satisfying the same inequality above.

Let $\omega $ be a basic form supported on the complement of a $2\rho $
-neighborhood of $\mathrm{Crit}\left( V\right) $, and let 
\begin{equation*}
\omega _{t}=\cos \left( tD_{B,s}\right) \omega =\frac{1}{2}\left(
e^{itD_{B,s}}+e^{-itD_{B,s}}\right) \omega ,
\end{equation*}
which is a solution to the generalized wave equation (\ref{basicwave})
corresponding to the operator $\Delta _{s}^{\prime }$ with initial
conditions $\omega _{0}=\omega $, $\frac{\partial }{\partial t}\omega _{0}=0$
. The family of forms $\omega _{t}$ is the unique solution to this
generalized wave equation as well, by the statements before and after (\ref
{basicwave}). Note that the formula above implies that $\omega _{t}$ is
basic.

By the unit propagation speed property of the basic wave equation
(Proposition \ref{UnitPropSpeed}), $\omega _{t}$ is identically zero on the $
\rho $-neighborhood of $\mathrm{Crit}\left( V\right) $ if $\left| t\right|
<\rho $. This implies that $D_{B,s}^{2}\omega _{t}=A\omega _{t}$ for $\left|
t\right| <\rho $, so that $\omega _{t}$ is the unique solution to the system 
\begin{equation*}
\frac{\partial ^{2}}{\partial t^{2}}\omega _{t}+A\omega
_{t}=0;\,\,\,\,\omega _{0}=\omega ,\frac{\partial }{\partial t}\omega _{0}=0.
\end{equation*}
We may therefore write $\omega _{t}=\cos \left( t\sqrt{A}\right) \omega $.

Let $\phi$ be a real-valued function with the following properties:

\begin{enumerate}
\item  $\phi $ is a positive even Schwarz function;

\item  $\phi \left( 0\right) =1$;

\item  the Fourier transform $\hat{\phi}\left( t\right) $ is supported in
the interval $\left[ -\rho ,\rho \right] $.
\end{enumerate}

For each nonnegative integer $m$, define $\phi_m$ by the formula 
\begin{equation*}
\phi_m(\lambda) = (1 + \lambda^2)^{2m}\phi(\lambda);
\end{equation*}
note that each $\phi_m$ satisfies (1) - (3).

For a basic form $\omega$ that is supported on the complement of the $2\rho$
-neighborhood of Crit$(V)$, 
\begin{eqnarray}
\phi_m \left( D_{B,s}\right) \omega &=&\frac{1}{2\pi }\int_{-\rho }^{\rho }%
\hat{\phi}_m\left( t\right) \left( e^{itD_{B,s}}\omega \right) \,dt  \notag
\\
&=&\frac{1}{2\pi }\int_{-\rho }^{\rho }\hat{\phi}_m\left( t\right) \left(
\cos \left( tD_{B,s}\right) \omega \right) \,dt\,\,\text{since }\hat{\phi}_m%
\text{ is even}  \notag \\
&=&\frac{1}{2\pi }\int_{-\rho }^{\rho }\hat{\phi}_m\left( t\right) \left(
\cos \left( t\sqrt{A}\right) \omega \right) \,dt  \notag \\
&=&\phi_m \left( \sqrt{A}\right) \omega .  \label{D_B,s and A}
\end{eqnarray}

The operator $\sqrt{A}$ is positive and has operator norm is bounded below
by $\frac{Cs}{\sqrt{2}}$ for $s$ sufficiently large. Thus, the operator norm
of $\phi_m \left( \sqrt{A}\right) $ (as an operator from $\mathcal{H}$ to
itself) is bounded above by 
\begin{equation*}
c_m\left( s\right) =\sup \left\{ \left| \phi_m \left( \lambda \right)
\right| \,:\,\lambda \geq \frac{Cs}{\sqrt{2}}\right\} .
\end{equation*}
It is clear that $c_m\left( s\right) $ is rapidly decreasing as $
s\rightarrow \infty $. By (\ref{D_B,s and A}), 
\begin{equation*}
\left\| \phi_m \left( D_{B,s}\right) \omega \right\| _{2}\leq c_m\left(
s\right) \left\| \omega \right\| _{2}
\end{equation*}
for every basic form $\omega $ supported on the complement of a $2\rho $
-neighborhood of $\mathrm{Crit}\left( V\right) $.

Next, let $L^{p}=L^{p}\left( \Omega _{B}^{*}\left( M\right) \right) $ denote
the $L^{p}-$norm closure of the space of smooth basic forms, and let $
W^{k}=W^{k}\left( \Omega _{B}^{*}\left( M\right) \right) $ denote the
closure of the space of such basic forms under the Sobolev $\left(
k,2\right) $-norm. By the basic elliptic estimates (Lemma \ref{elliptic2}),
the operator norm of $\left( 1+D_{B,s}^{2}\right) ^{-1}:W^{k}\longrightarrow
W^{k+2}$ is bounded by a polynomial in $s$. The basic version of the Sobolev
imbedding theorem (Lemma \ref{sobolev}) implies that $\left(
1+D_{B,s}^{2}\right) ^{-k}:L^{2}\longrightarrow L^{\infty }$ is a bounded
map whose operator norm is bounded by a polynomial in $s$ if $k>\frac{n}{4}$
. Using basic duality (Lemma \ref{dualL^1}) and essential self-adjointness
of $D_{B,s}$, we see that $\left( 1+D_{B,s}^{2}\right)
^{-k}:L^{1}\longrightarrow L^{2}$ is also a bounded map whose operator norm
is bounded by a polynomial in $s$ whenever $k>\frac{n}{4}$. Note that all of
the statements above hold for the operator $A$ as well as for $D_{B,s}^{2}$.
Now, given a basic form $\omega $ supported on the complement of a $2\rho $
-neighborhood of $\mathrm{Crit}\left( V\right) $ and $k>\frac{n}{4}$, 
\begin{multline*}
\left\| \phi \left( D_{B,s}\right) \omega \right\| _{\infty }=\left\| \phi
\left( \sqrt{A}\right) \omega \right\| _{\infty } \\
\left\| \left( 1+A\right) ^{-k}\left( 1+A\right) ^{k}\phi \left( \sqrt{A}%
\right) \left( 1+A\right) ^{k}\left( 1+A\right) ^{-k}\omega \right\|
_{\infty } \\
\leq \left\| \left( 1+A\right) ^{-k}\right\| _{L^{2}\rightarrow L^{\infty
}}c_{k}\left( s\right) \left\| \left( 1+A\right) ^{-k}\omega \right\| _{2} \\
\leq \left\| \left( 1+A\right) ^{-k}\right\| _{L^{2}\rightarrow L^{\infty
}}c_{k}\left( s\right) \left\| \left( 1+A\right) ^{-k}\right\|
_{L^{2}\rightarrow L^{1}}\left\| \omega \right\| _{1} \\
\leq p\left( s\right) c_{k}\left( s\right) \left\| \omega \right\| _{1}
\end{multline*}
where $p\left( s\right) $ is a polynomial in $s$. Next, since $\phi $ is
rapidly decreasing, $\phi \left( D_{B,s}\right) $ has a continuous basic
kernel $k_{B,s}\left( x,y\right) $ (Lemma \ref{smoothkernel}), and we have
the inequality 
\begin{equation*}
\left\| k_{B,s}\left( x,\cdot \right) \right\| _{\infty }\leq \sup_{\int
\left\| \omega \right\| =1,x\in M}\left\| \int_{M}k_{B,s}\left( x,y\right)
\,\omega \left( y\right) \omega _{\mathrm{vol}}\left( y\right) \right\|
\,\leq p\left( s\right) c_{k}\left( s\right)
\end{equation*}
from the above. Thus, $k_{B,s}\left( x,y\right) \rightarrow 0$ uniformly as $
s\rightarrow \infty $.
\end{proof}

Let $\left( M,\mathcal{F}\right) $ be a Riemannian foliation, let $V$ be a $ 
\mathcal{F}$-nondegenerate basic vector field, let $L$ be a leaf closure of $
\mathcal{F}$, and let $\mathcal{O}_{L}=$ $\mathcal{O}_{L}\left( V\right) $
denote the orientation line bundle of $V$ at $L$. We denote by $\Omega
^{*}\left( L,\mathcal{F},\mathcal{O}_{L}\right) $ the space of differential
forms on $L$ with values in $\mathcal{O}_{L}$; there is a well-defined
differential $d$ on this space that is simply the exterior derivative when $ 
\mathcal{O}_{L}$ is trivial (see \cite[Section I.7]{BT}). Locally an element
of $\Omega ^{*}\left( L,\mathcal{F},\mathcal{O}_{L}\right) $ can be written
as a sum of tensors $\omega \otimes s$, where $\omega $ is an ordinary
differential form on $M$ and $s$ is a smooth section of $\mathcal{O}_{L}$.
Let $X$ be a vector field on $L$. Then we can locally define interior
multiplication $i(X)$ on $\Omega ^{*}\left( L,\mathcal{F},\mathcal{O}
_{L}\right) $ by decreeing that $i(X)(\omega \otimes s)=(i(X)\omega )\otimes
s$ and extending linearly. It is straightforward to check that this
definition of interior product is independent of the trivialization of $ 
\mathcal{O}_{L}\left( V\right) $, and so we may therefore make the following
definitions:

\begin{definition}
\label{twisted}Let $M$, $\mathcal{F}$, $V$, and $L$ have the properties
listed in the previous paragraph. The \textrm{space of basic differential
forms with values in $\mathcal{O}_{L}$} is denoted $\Omega _{B}^{\ast
}\left( L,\mathcal{F},\mathcal{O}_{L}\right) $, and is defined to be the
subcomplex of forms $\alpha $ in $\Omega ^{\ast }\left( L,\mathcal{F}, 
\mathcal{O}_{L}\right) $ with the property that $i(X)\alpha =0$ and $
i(X)d\alpha =0$ for every vector field $X$ on $L$ tangent to the leaves of $ 
\mathcal{F}$ restricted to $L$. The cohomology of this subcomplex is denoted 
$H_{B}^{\ast }\left( L,\mathcal{F},\mathcal{O}_{L}\right) $ and called the 
\textrm{basic deRham cohomology of $L$ with values in $\mathcal{O}_{L}$}.
Finally, the \textrm{basic Euler characteristic of $L$ with values in $%
\mathcal{O}_{L}$} is defined by the formula $\chi _{B}\left( L,\mathcal{F}, 
\mathcal{O}_{L}\right) =\sum_{k}(-1)^{k}\mathrm{dim\,}H_{B}^{k}\left( L, 
\mathcal{F},\mathcal{O}_{L}\right) .$
\end{definition}

\begin{theorem}
\label{BasicHopf}(Basic Hopf Index Theorem) Let $M$, $\mathcal{F}$, and $V$
be as in Definition \ref{twisted}. For each critical leaf closure $L$, let $ 
\mathrm{ind}_{L}\left( V\right) $ be the index of $V$ at $L$. Then 
\begin{equation*}
\chi _{B}\left( M,\mathcal{F}\right) =\sum_{L\,\mathrm{critical} }
  \mathrm{%
ind}_{L}\left( V\right)\chi _{B}\left( L,\mathcal{%
F},\mathcal{O}_{L}\right) .
\end{equation*}
\end{theorem}

\begin{proof}
Without loss of generality, we may assume that $V$ is a basic normal vector
field (otherwise, project to $N\mathcal{F}$). By Proposition~\ref{indepofs},
we have that $\chi _{B}\left( M,\mathcal{F}\right) =\lim_{s\rightarrow
\infty }\left( \mu _{s}^{\mathrm{even}}-\mu _{s}^{\mathrm{odd}}\right) $,
which may be obtained by integrating the traces of the corresponding kernels
of $\left. \phi \left( D_{B,s}^{2}\right) \right| _{L^{2}\left( \Omega _{B}^{%
\mathrm{even}}\left( M\right) \right) }$ and $\left. \phi \left(
D_{B,s}^{2}\right) \right| _{L^{2}\left( \Omega _{B}^{\mathrm{odd}}\left(
M\right) \right) }$. By Proposition~\ref{complement}, the kernels of these
operators go to zero uniformly on the complement of a fixed but arbitrarily
small neighborhood of the critical leaf closures of $V$. For each critical
leaf closure $L$, let $\psi _{L}$ be a smooth, radial, basic function that
is identically $1$ in a tubular neighborhood of radius $2\rho $ around $L$
and supported within a tubular neighborhood of radius $3\rho $ (assume that
we have chosen $\rho $ small enough so that this is possible for each $L$).
Then we have that 
\begin{equation*}
\chi _{B}\left( M,\mathcal{F}\right) =\sum_{L\,\mathrm{critical}%
}\lim_{s\rightarrow \infty }\mathrm{tr}\left( \left. \psi _{L}\phi \left(
D_{B,s}^{2}\right) \right| _{L^{2}\left( \Omega _{B}^{\mathrm{even}}\left(
M\right) \right) }\right) -\mathrm{tr}\left( \left. \psi _{L}\phi \left(
D_{B,s}^{2}\right) \right| _{L^{2}\left( \Omega _{B}^{\mathrm{odd}}\left(
M\right) \right) }\right) .
\end{equation*}
We now use Lemma~\ref{deform} to observe that $V$ may be deformed without
changing $\mathrm{ind}_{L}\left( V\right) $ or $\mathcal{O}_{L}$ so that
within a tubular neighborhood of radius $4\rho $ around $L$, $V=\nabla f$
for a basic function $f$, such that if $\mathrm{ind}_{L}\left( V\right) =+1$
then $f$ has even Morse index and if $\mathrm{ind}_{L}\left( V\right) =-1$
then $f$ has odd Morse index; again we possibly decrease $\rho $ so that the
conclusion of this proposition holds. A unit propagation speed argument
shows that the traces are independent of the choice of $V$ with those
properties, and thus we may calculate the contributions from each tubular
neighborhood as if $V=\nabla f$. By the results of \cite{A}, if $D_{L}^{2}$
is the basic Laplacian on $L$ with coefficients in $\mathcal{O}_{L}$ , then 
\begin{equation*}
\lim_{s\rightarrow \infty }\mathrm{tr}\left( \left. \psi _{L}\phi \left(
D_{B,s}^{2}\right) \right| _{L^{2}\left( \Omega _{B}^{\mathrm{even}}\left(
M\right) \right) }\right) =\mathrm{tr}\left( \left. \phi \left(
D_{L}^{2}\right) \right| _{L^{2}\left( \Omega _{B}^{\mathrm{even}+\frac{%
\mathrm{ind}_{L}\left( V\right) -1}{2}}\left( L,\mathcal{F},\mathcal{O}%
_{L}\right) \right) }\right) ,
\end{equation*}
with the analogous result for the odd case. Thus, 
\begin{eqnarray*}
\chi _{B}\left( M,\mathcal{F}\right) &=&\sum_{L\,\mathrm{critical}}\mathrm{tr%
}\left( \left. \phi \left( D_{L}^{2}\right) \right| _{L^{2}\left( \Omega
_{B}^{\mathrm{even}+\frac{\mathrm{ind}_{L}\left( V\right) -1}{2}}\left( L,%
\mathcal{F},\mathcal{O}_{L}\right) \right) }\right) - \\
&&\mathrm{tr}\left( \left. \phi \left( D_{L}^{2}\right) \right|
_{L^{2}\left( \Omega _{B}^{\mathrm{odd}+\frac{\mathrm{ind}_{L}\left(
V\right) -1}{2}}\left( L,\mathcal{F},\mathcal{O}_{L}\right) \right) }\right)
\\
&=&\sum_{L\,\mathrm{critical}}\dim H_{B}^{\mathrm{even}+\frac{\mathrm{ind}
_{L}\left( V\right) -1}{2}}\left( L,\mathcal{F},\mathcal{O}_{L}\right) -\dim
H_{B}^{\mathrm{odd}+\frac{\mathrm{ind}_{L}\left( V\right) -1}{2}}\left( L,%
\mathcal{F},\mathcal{O}_{L}\right) \\
&=&\sum\begin{Sb} L\,\mathrm{critical}  \\ \mathrm{ind}_{L}\left( V\right)
=+1  \end{Sb}  \chi _{B}\left( L,\mathcal{F},\mathcal{O}_{L}\right) -\sum 
\begin{Sb} L\,\mathrm{critical}  \\ \mathrm{ind}_{L}\left( V\right) =-1 
\end{Sb}  \chi _{B}\left( L,\mathcal{F},\mathcal{O}_{L}\right) .
\end{eqnarray*}
\end{proof}

\begin{remark}
If all of the bundles $\mathcal{O}_{L}$ are trivial and the critical leaf
closures are leaves, then the formula simplifies to 
\begin{equation*}
\chi _{B}\left( M,\mathcal{F}\right) =\sum_{L\,\mathrm{critical}}\mathrm{ind}
_{L}\left( V\right) ,
\end{equation*}
a formula which has the precise form of the ordinary Hopf Index Theorem.
\end{remark}

\begin{corollary}
\label{obstruction} Suppose that there exists a basic vector field on a
Riemannian foliation $(M,\mathcal{F})$ that is nowhere tangent to $\mathcal{F%
}$. Then $\chi _{B}\left( M,\mathcal{F}\right) =0.$
\end{corollary}\qed

\begin{corollary}
\label{fundgroup} Suppose that there exists a basic vector field on a
Riemannian foliation $(M,\mathcal{F})$ that is nowhere tangent to $\mathcal{F%
}$ and that the codimension of $\mathcal{F}$ is less than $3$. Then $M$ has
infinite fundamental group.
\end{corollary}

\begin{proof}
The hypotheses imply 
\begin{eqnarray*}
\chi _{B}\left( M,\mathcal{F}\right) &=&\dim
H_{B}^{0}\left( M,\mathcal{F}\right) -\dim H_{B}^{1}\left( M,\mathcal{F}
\right) +\dim H_{B}^{2}\left( M,\mathcal{F}\right) \\
&=&0.
\end{eqnarray*}
 Since $\dim
H_{B}^{0}\left( M,\mathcal{F}\right) \geq 1$, $\dim H_{B}^{1}\left( M, 
\mathcal{F}\right) \geq 1$. Since $H_{B}^{1}\left( M,\mathcal{F}\right) $
injects into $H^{1}(M)$ \cite[Proposition 4.1]{To}, Poincar\'{e} duality and
the Universal Coefficient Theorem imply that $\mathrm{rank}\,H_{1}(M)\geq 1$
. The result follows.
\end{proof}

We now illustrate Theorem \ref{BasicHopf} in the following examples. In each
of these examples, there does not exist a basic vector field that is nowhere
tangent to the foliation.

\begin{example}
Consider the sphere $S^{2}$ in spherical coordinates $\left( \theta ,\varphi
\right) \in \left[ 0,2\pi \right] \times \left[ 0,\pi \right] $. Let $\alpha 
$ be a fixed irrational multiple of $\pi $, and let $M=\Bbb{R\times S}
^{2}
\slash
 \sim $, where $\left( t,\theta ,\varphi \right) \sim \left(
t+1,\theta +\alpha ,\varphi \right) $. The $t$-parameter curves make $M$
into a codimension-2 foliation $\mathcal{F}$, which is Riemannian with
respect to the standard metric. Observe that the leaf closures are level
sets where $\varphi $ is constant. Consider the vector field 
\begin{equation*}
V=\left( \cos \varphi \right) \partial _{\theta }+\left( \sin \varphi \cos
\varphi \right) \partial _{\varphi }.
\end{equation*}
This vector field is invariant under rotations in $\theta $ and is smooth on 
$S^{2}$, since it is the restriction of 
\begin{equation*}
\left( xz^{2}-yz\right) \partial _{x}+\left( yz^{2}+xz\right) \partial
_{y}+\left( -x^{2}z-y^{2}z\right) \partial _{z}
\end{equation*}
on $\Bbb{R}^{3};$ therefore, it is a smooth basic vector field on $\left( M, 
\mathcal{F}\right) $. The critical leaf closures for this vector field
correspond to the poles $\left( z=\pm 1\right) $ and the equator $\left(
\varphi =\frac{\pi }{2}\right) $. At the north pole leaf $L_{1}$ ($z=1$),
the matrix for $V_{L_{1}}$ is $\left( 
\begin{array}{ll}
1 & -1 \\ 
1 & 1
\end{array}
\right) $ in $\left( 
\begin{array}{l}
x \\ 
y
\end{array}
\right) $ coordinates, so that $\mathrm{ind}_{L_{1}}\left( V\right) =1$. At
the south pole $L_{-1}$ ($z=-1$), the matrix for $V_{L_{-1}}$ is $\left( 
\begin{array}{ll}
-1 & 1 \\ 
-1 & -1
\end{array}
\right) $ in $\left( 
\begin{array}{l}
x \\ 
y
\end{array}
\right) $ coordinates, so that $\mathrm{ind}_{L_{-1}}\left( V\right) =1$.
For these leaf closures, we have $\chi _{B}\left( L_{\pm 1},\mathcal{F}, 
\mathcal{O}_{L_{\pm 1}}\right) =1$ since $\mathcal{O}_{L_{\pm 1}}$ and $ 
\mathcal{F}$ are trivial. At the equator $L_{0}$, $V_{L_{0}}$ is
multiplication by $-1$ on each $1$-dimensional normal space to $L_{0}$, so
that $\mathrm{ind}_{L_{0}}\left( V\right) =-1$. The orientation bundle is
the trivial co-normal bundle. Next, observe that this leaf closure is a flat
torus, on which the foliation restricts to be the irrational flow. Since the
vector field $\partial _{\theta }$ is basic, nonsingular, and orthogonal to
the foliation on this torus, 
\begin{equation*}
\chi _{B}\left( L_{0},\mathcal{F},\mathcal{O}_{L_{0}}\right) =\chi
_{B}\left( L_{0},\mathcal{F}\right) =0
\end{equation*}
by Corollary \ref{obstruction}. By Theorem \ref{BasicHopf}, we conclude that 
\begin{eqnarray*}
\chi _{B}\left( M,\mathcal{F}\right) &=&\sum\begin{Sb} L\,\mathrm{critical} 
\\ \mathrm{ind}_{L}\left( V\right) =+1  \end{Sb}  \chi _{B}\left( L,\mathcal{%
F},\mathcal{O}_{L}\right) -\sum\begin{Sb} L\,\mathrm{critical}  \\ \mathrm{%
ind} _{L}\left( V\right) =-1  \end{Sb}  \chi _{B}\left( L,\mathcal{F},%
\mathcal{O}_{L}\right) \\
&=&\left( 1+1\right) -\left( 0\right) =2.
\end{eqnarray*}

In this example, one may independently verify that $\dim H_{B}^{k}\left( M, 
\mathcal{F}\right) =1$ for $k=0$ or $k=2$, since $M$ is a transversally
oriented, taut, codimension-2 Riemannian foliation (see \cite{To}). Also,
every closed basic one-form can be written as $g\left( \varphi \right)
d\varphi $ for a smooth function $g$ such that $\frac{\partial ^{m}g}{%
\partial \varphi ^{m}}=0$ for even $m$ at $\varphi =0$ or $\pi $. Every
exact basic one-form is the differential of a smooth function $h$ such that $
\frac{\partial ^{m}h}{\partial \varphi ^{m}}=0$ for odd $m$ at $\varphi =0$
or $\pi $; thus, $\dim H_{B}^{1}\left( M,\mathcal{F}\right) =0.$ We have
therefore directly computed that $\chi _{B}\left( M,\mathcal{F}\right)
=1-0+1=2.$ Observe that Corollary \ref{obstruction} implies that there does
not exist a basic vector field on $\left( M,\mathcal{F}\right) $ that is
nowhere tangent to the leaves.
\end{example}

\begin{example}
Let $X$ be any smooth, closed manifold whose fundamental group is the free
product $\Bbb{Z}*\Bbb{Z}$ (such as the connected sum of two copies of $
S^{2}\times S^{1}$); choose a metric so that the volume of $X$ is one. Let $ 
\widetilde{X}$ denote the universal cover of $X$. Suspend an action of the
fundamental group on the torus to get the manifold $Y=\widetilde{X}\times
S^{1}\times S^{1}
\slash
 \sim $, where the equivalence relation is defined
as follows. Let $a$ and $b$ be the generators of $\pi _{1}\left( X\right) ,$
which act on the right by isometries on $\widetilde{X}$. Let $\alpha $ be an
irrational multiple of $2\pi $, and let 
\begin{eqnarray*}
&&\left( x,\theta _{1},\theta _{2}\right) \sim \left( xa,\theta _{1}+\alpha
,\theta _{2}\right) ,\text{ and } \\
&&\left( x,\theta _{1},\theta _{2}\right) \sim \left( xb,2\pi -\theta
_{1},2\pi -\theta _{2}\right) .
\end{eqnarray*}
These actions on the torus (with the flat metric) are orientation-preserving
isometries, so that the $x$-parameter (immersed) submanifolds form a
transversally oriented, Riemannian foliation $\left( Y,\mathcal{F}\right) $.
The leaf closures are parametrized by $\theta _{2}\in \left[ 0,\pi \right] $
; observe that the leaf closures corresponding to $\theta _{2}=0$ and $
\theta _{2}=\pi $ are not transversally oriented, even though $\mathcal{F}$
is.

We first calculate $\chi _{B}\left( Y,\mathcal{F}\right) $ directly. As in
the previous example, we observe that $\dim H_{B}^{k}\left( Y,\mathcal{F}
\right) =1$ for $k=0$ or $k=2$, since $M$ is a transversally oriented, taut,
codimension-2 Riemannian foliation. Each closed, basic one-form can be
written as $g\left( \theta _{2}\right) \,d\theta _{2}$, where $g$ is a
smooth function on $S^{1}$such that $g\left( 2\pi -\theta _{2}\right)
=-g\left( \theta _{2}\right) $. Since this is the differential of the
well-defined basic function $f\left( \theta _{2}\right) =\int_{0}^{\theta
_{2}}g\left( t\right) \,dt$, we see that $\dim H_{B}^{1}\left( Y,\mathcal{F}
\right) =0$. Therefore, $\chi _{B}\left( Y,\mathcal{F}\right) =1-0+1=2.$

Next, consider the basic vector field $W=\sin \left( \theta _{2}\right)
\partial _{\theta _{2}}$. This vector field is singular at the leaf closures 
$\theta _{2}=0$ and $\theta _{2}=\pi $ (both of codimension 1), and it is
orthogonal to the leaves everywhere. The index of $W$ at $\theta _{2}=0$ is $
1$, and its index at $\theta _{2}=\pi $ is $-1$. The polar decomposition of
the linearization of $W$ at $\theta _{2}=0$ is simply $1\ast 1,$ so there
are no negative eigenvectors of the orthogonal part. Thus, the orientation
line bundle is trivial at $\theta _{2}=0.$ Therefore, 
\begin{equation*}
\chi _{B}\left( \left\{ \theta _{2}=0\right\} ,\mathcal{F},\mathcal{O}%
_{L}\right) =\chi _{B}\left( \left\{ \theta _{2}=0\right\} ,\mathcal{F}%
\right) =1-0=1;
\end{equation*}
note that there are no closed one-forms on this leaf closure, because it is
not transversally oriented. The polar decomposition of the linearization of $
W$ at $\theta _{2}=\pi $ is $1\ast \left( -1\right) $, and the orientation
line bundle $\mathcal{O}_{\pi }$ is simply the normal bundle to the leaf
closure. This bundle is nontrivial and has no basic sections, so $\dim
H_{B}^{0}\left( \left\{ \theta _{2}=0\right\} ,\mathcal{F},\mathcal{O}
_{L}\right) =0$. On the other hand, let $s$ denote a basic section of the
pullback of $\mathcal{O}_{\pi }$ via the lift of the leaf closure $\theta
_{2}=\pi $ to $\widetilde{X}\times S^{1}\times S^{1}$; such a section exists
because the pullback of $\mathcal{O}_{\pi }$ is trivial. The basic one-form $
d\theta _{1}\otimes s$ is closed and descends to a basic one-form on $Y$
with values in $\mathcal{O}_{\pi }$, because the orientation-reversing
action of $b$ changes the sign of both $d\theta _{1}$ and $s$. It follows
that $\dim H_{B}^{1}\left( \left\{ \theta _{2}=0\right\} ,\mathcal{F}, 
\mathcal{O}_{L}\right) =1$. Thus, 
\begin{equation*}
\chi _{B}\left( \left\{ \theta _{2}=0\right\} ,\mathcal{F},\mathcal{O}%
_{L}\right) =0-1=-1.
\end{equation*}
By Theorem \ref{BasicHopf}, 
\begin{eqnarray*}
\chi _{B}\left( Y,\mathcal{F}\right) &=&\sum\begin{Sb} L\,\mathrm{critical} 
\\ \mathrm{ind}_{L}\left( V\right) =+1  \end{Sb}  \chi _{B}\left( L,\mathcal{%
F},\mathcal{O}_{L}\right) -\sum\begin{Sb} L\,\mathrm{critical}  \\ \mathrm{%
ind}_{L}\left( V\right) =-1  \end{Sb}  \chi _{B}\left( L,\mathcal{F},%
\mathcal{O}_{L}\right) \\
&=&\left( 1\right) -\left( -1\right) =2.
\end{eqnarray*}
These calculations show that there does not exist a basic vector field on $
\left( Y,\mathcal{F}\right) $ that is nowhere tangent to the foliation; in
fact, any $\mathcal{F}$-nondegenerate basic vector field on $\left( Y, 
\mathcal{F}\right) $ must be tangent to the foliation on at least two
distinct leaf closures.
\end{example}

\begin{example}
Let $N$ be a smooth, four-dimensional closed manifold with finite
fundamental group. Let $\left( M,\mathcal{F}\right) $ be any foliation that
is obtained by suspending a discrete subgroup $\Gamma $ of a compact Lie
group of diffeomorphisms of $N$. That is, choose a manifold $X$ along with a
surjective homomorphism $\phi :\pi _{1}\left( X\right) \rightarrow \Gamma $,
and let $M=\widetilde{X}\times N
\slash
 \pi _{1}\left( X\right) $, where $
\pi _{1}\left( X\right) $ acts on the universal cover $\widetilde{X}$ by
deck transformations and on $N$ via $\phi $. The foliation $\mathcal{F}$ is
locally given by the $\widetilde{X}$-parameter submanifolds. \newline
Choose metrics for $X$ and $N$; by averaging the metric on $N$ over the Lie
group of diffeomorphisms, we may and do assume that $\pi _{1}\left( X\right) 
$ acts on $N$ by isometries. The metric on $M$ defined locally as the
product of these metrics is bundle-like for the foliation $\mathcal{F}$.
Furthermore, this foliation is taut, so the standard form of Poincar\'{e}
duality holds for basic cohomology (see \cite{To}). The basic forms of $
\left( M,\mathcal{F}\right) $ are given by forms on $N$ that are invariant
under the discrete group of isometries, and the basic cohomology is
isomorphic to the cohomology of invariant forms on $N$. Since $\pi
_{1}\left( N\right) $ is finite, $H_{B}^{1}\left( M,\mathcal{F}\right) $ is
trivial. We conclude that 
\begin{eqnarray*}
\chi _{B}\left( M,\mathcal{F}\right) &=&\sum_{j=0}^{4}\left( -1\right)
^{j}\dim H_{B}^{j}\left( M,\mathcal{F}\right) \\
&=&\dim H_{B}^{0}\left( M,\mathcal{F}\right) +\dim H_{B}^{2}\left( M,%
\mathcal{F}\right) +\dim H_{B}^{4}\left( M,\mathcal{F}\right) \\
&=&2+\dim H_{B}^{2}\left( M,\mathcal{F}\right) \geq 2.
\end{eqnarray*}
Theorem \ref{BasicHopf} implies that every $\mathcal{F}$-nondegenerate basic
vector field on $\left( M,\mathcal{F}\right) $ must have at least two
distinct critical leaf closures.\newline
\end{example}


\begin{thebibliography}{99}
\bibitem{A}  J. A. Alvarez Lop\'{e}z, \emph{Morse Inequalities for
Pseudogroups of Local Isometries}, J. Diff. Geom. \textbf{37}(1993), 603-638.

\bibitem{BT}  R. Bott, L. W. Tu, \emph{Differential Forms in Algebraic
Topology}, New York: Springer-Verlag, 1982.

\bibitem{D}  J. Dieudonn\'{e}, \emph{Treatise on Analysis VIII} (trans. by
L. Fainsilber), Boston: Academic Press, 1993.

\bibitem{EKHS}  A. El Kacimi, G. Hector, V. Sergiescu, \emph{La\
cohomologie\ basique\ d'un\ feuilletage\ riemannien\ est \ de \ dimension\
finie}, Math. Z. \textbf{188}(1985), 593--599.

\bibitem{ESW}  A. El Soufi, X. P. Wang, \emph{Some\ remarks\ on\ Witten's\
method.\ Poincar\'{e}-Hopf\ theorem\ and\ Atiyah-Bott\ formula}, Ann. Global
Anal. Geom. \textbf{5}(1987), 161-178.

\bibitem{Gh}  E. Ghys, \emph{Un feuilletage analytique dont la cohomologie
basique est de dimension infinie}, Publ. de l'IRMA, Lille \textbf{7}(1985).

\bibitem{KT}  F. W. Kamber, and Ph. Tondeur, \emph{de Rham-Hodge theory for
Riemannian foliations,} Math. Ann. \textbf{277}(1987), 415--431.

\bibitem{La}  O. A. Ladyzhenskaya, \emph{The Boundary Value Problems of
Mathematical Physics} (trans. by J. Lohwater), New York: Springer-Verlag,
1985.

\bibitem{LeR}  J. M. Lee, K. Richardson, \emph{Riemannian foliations and
eigenvalue comparison,} Ann. Glob. Anal. Geom. \textbf{16}(1998), 497--525.

\bibitem{Ma}  X. M. Masa, \emph{Duality and minimality in Riemannian
foliations}, Comment. Math. Helv. \textbf{67}(1992), 17--27.

\bibitem{Mi}  S. Mizohata, \emph{The theory of partial differential equations%
}, London: Cambridge University Press, 1973.

\bibitem{Mo}  P. Molino, \emph{Riemannian foliations}, Boston: Birkh\"{a}
user, 1988.

\bibitem{PR}  E. Park, K. Richardson, \emph{The basic Laplacian of a
Riemannian foliation,} Amer. J. Math. \textbf{118}(1996), 1249--1275.

\bibitem{R}  K. Richardson, \emph{Traces of heat operators on Riemannian
foliations}, preprint.

\bibitem{Ro}  J. Roe, \emph{Elliptic operators, topology, and asymptotic
methods}, Pitman Research Notes in Math. \textbf{179}, Longman Scientific
and Technical: Harlow, 1988.

\bibitem{Ta}  M. E. Taylor, \emph{Partial Differential Equations I}, New
York: Springer-Verlag, 1996.

\bibitem{To}  Ph. Tondeur, \emph{Geometry of Foliations}, Basel: Birkh\"{a}
user, 1997.

\bibitem{W}  E. Witten, \emph{Supersymmetry and Morse theory}, J. Diff.
Geom. \textbf{17}(1982), 661-692.

\bibitem{Z}  A. Zeggar, \emph{Nombre de Lefschetz basique pour un feulletage
riemannien}, Ann. Fac. Sci. Toulouse Math., \textbf{1}(1992), 105--131.
\end{thebibliography}
\end{document}